\documentstyle[12pt]{article}

\makeatletter

\@addtoreset{equation}{section}
\def\section{\@startsection {section}{1}{\z@}{-3.5ex plus -1ex minus 
 -.2ex}{2.3ex plus .2ex}{\normalsize\bf}}
\def\@maketitle{\newpage
 \null
 \vskip 2em \begin{center}
 {\large \@title \par} \vskip 1.5em {\normalsize \lineskip .5em
\begin{tabular}[t]{c}\@author 
 \end{tabular}\par} 
 \vskip 1em {\normalsize \@date} \end{center}
 \par
 \vskip 1.5em} 
\makeatother
\newtheorem{th}{Theorem}[section]
\newtheorem{lm}[th]{Lemma}
\newtheorem{prop}[th]{Proposition}
\newtheorem{df}[th]{Definition}

\newtheorem{cl}[th]{Claim}
\newtheorem{cor}[th]{Corollary}

\def\tr{{\rm tr}}
\def\e{\varepsilon}

\input amssym.def
\newsymbol\rtimes 226F
\newsymbol\varnothing 203F
\def\emptyset{\varnothing} 
\newsymbol\varsubsetneq 2020 
\begin{document}

\title{\bf Central sequence subfactors and double commutant properties}
\author{Keiko Kawamuro\\
        Department of Mathematical Sciences \\
        University of Tokyo, Komaba, Tokyo, 153, JAPAN  \\
        e-mail: yuri@ms.u-tokyo.ac.jp }
%\footnotetext{1991 {\it Mathematics Subject Classification} 46L37.}
%\footnotetext[2]{This manuscript is written by using \LaTeX.}
\date{}
\maketitle
\date{}
\begin{abstract}

First, we construct the Jones tower and tunnel of
the central sequence subfactor
arising from a hyperfinite type II$_1$ subfactor with
finite index and  finite depth,
and prove each algebra has the double commutant property
in the ultraproduct of the enveloping II$_1$ factor.
Next, we show the equivalence between 
Popa's strong amenability and the double commutant property
of the central sequence factor for subfactors as above without
assuming the finite depth condition.
\end{abstract}

\section{Introduction}

Subfactor theory was initiated by V. F. R. Jones in the 80's (\cite{J2}).
The central sequence subfactor, which is the key concept in this paper,
is one of the 
important notions in the analytical approaches to subfactor theory.
We review the historical background.

Let $R$ be the hyperfinite II$_1$ factor, $G$ a finite group
and $\alpha$ an outer action of $G$ on $R.$
We fix a free ultrafilter $\omega$ on ${\bf N}.$
Originally, it has been known that 
the classification of group actions can be understood conceptually
in terms of the subfactor $R_{\omega}^G \subset R_{\omega},$
%$B$D$^$j!$(Bcentral sequence algebra$B>e$G!$(BRohlin
%$B$J$I$r9M$($k$3$H$,!$(BR_{\omega}^G $B$H(B R_{\omega}$B$N4X78$rD4$Y$F$$$k!$$H(B
%$B9M$($k$3$H$,$G$-$k!$$H8@$&DxEY$N0UL#$G$9!%(B
%
%the classification of group actions is related to studying the subfactor
%$R_{\omega}^G \subset R_{\omega},$
where $R_{\omega}$ stands for a central sequence algebra 
$R^{\omega} \cap R'.$
We considered the following subfactor of S. Popa (\cite{Paction})
$$N := \left\{ 
\begin{array}{c|c}
\left(
\begin{array}{ccc}
\alpha_{g_1}(x) & & \\
  & \ddots & \\
 &  & \alpha_{g_n}(x)
\end{array}
\right) & x \in R 
\end{array}
\right\}
\subset M := R \otimes M_n({\bf C}),$$
where $\alpha_g$ 
is the image of $g \in G$ under a homomorphism (``action'') $\alpha$
of $G$ into the group ${\rm Aut}(R)$ of $*$-automorphism of $R.$
(This example was initially hinted by Jones.)
%means $$\alpha : G \owns g \mapsto \alpha_g \in {\rm Aut} (R).$$
Then we have
$N^{\omega} \cap M' = R_{\omega}^G$
and $M_{\omega} = R_{\omega}.$
This is a special case of the {\it central sequence subfactor}
$N^{\omega} \cap M' \subset M^{\omega} \cap M'.$
%The exact definition of it is as follows. \\
%{\bf Definition}  Let $N \subset M$ be a hyperfinite type II$_1$ subfactor
%with finite index and finite depth.
%We call the subfactor $N^{\omega} \cap M' \subset M^{\omega} \cap M'$
%the {\it central sequence subfactor}. \\
The central sequence subfactor has been introduced by A. Ocneanu 
with intention of generalizing the classification theory
of group actions on factors to that of ``actions'' of new
algebraic objects called {\it paragroups} (\cite{O1}, \cite{EK}).
%A. Ocneanu has introduced paragroups as generalized group.
%He has payed attention to the central sequence subfactor
%on purpose of classifying the paragroup actions.
That is,
let $N \subset M$ be a hyperfinite type II$_1$ subfactor
with finite index and finite depth, and
we construct a paragroup ${\cal G}$ from 
$N \subset M$, then the idea of the paragroup theory is to regard
$N \subset M$ as ``$ N \subset N \rtimes {\cal G}.$''

S. Popa has proved in \cite{Pacta} that we can reconstruct the original 
subfactor from the paragroup.  
This is the {\it generating property} of Popa.
It gives another proof of the uniqueness of 
outer actions of a finite group on the hyperfinite 
II$_1$ factor.
He has also defined a notion of amenability of subfactors
based on an analogy with groups and proved
the equivalence of strong amenability and the generating property
for hyperfinite type II$_1$ subfactors with finite indices.
This theorem of Popa also implies several classification results
about group actions on factors.

We obtain another type II$_1$ subfactor from $N \subset M,$
which is called the {\it asymptotic inclusion.}
It is defined by $M \vee (M' \cap M_{\infty}) \subset M_{\infty},$
where we denote the enveloping algebra of $N \subset M$
by $M_{\infty}.$
Recently it has been studied in various fields
such as topological quantum field theory (\cite[Section 12]{EK}), 
sector theory (\cite{LR}),
quantum doubles and so on.

Several properties 
of the central sequence subfactor 
have been studied by Y. Kawahigashi and Ocneanu.
The following theorem of Ocneanu is especially
important.
(\cite{O2}, \cite[Theorem 15.32]{EK}.)
We discuss it later in Section \ref{application}.
\medskip \\
{\bf Theorem (Ocneanu)}
The paragroups of the 
central sequence subfactor and of the asymptotic inclusion 
are mutually dual. 
\medskip \\
Looking at the proof of the theorem, we notice that
$( M^{\omega} \cap M')' \cap M_{\infty}^{\omega} \subset 
( N^{\omega} \cap M')' \cap M_{\infty}^{\omega}$
has the same higher relative commutant as the asymptotic inclusion,
and 
\begin{eqnarray*}
N^{\omega} \cap M'&=& (( N^{\omega} \cap M')' \cap
M_{\infty}^{\omega})' \cap M_{\infty}^{\omega},  \\
M^{\omega} \cap M' &=& (( M^{\omega} \cap M')' \cap
M_{\infty}^{\omega})' \cap M_{\infty}^{\omega}.
\end{eqnarray*}
We construct the Jones towers and tunnels from
$$N^{\omega} \cap M' \subset M^{\omega} \cap M'$$ and
$$( M^{\omega} \cap M')' \cap M_{\infty}^{\omega} \subset 
( N^{\omega} \cap M')' \cap M_{\infty}^{\omega}$$ 
in $M_{\infty}^{\omega}.$
We set 
$P_0 := N^{\omega} \cap M',$
$P_1 := M^{\omega} \cap M',$
$P_0^c := ( N^{\omega} \cap M')' \cap M_{\infty}^{\omega} $ and
$P_1^c := ( M^{\omega} \cap M')' \cap M_{\infty}^{\omega},$
where ``c'' stands for the relative commutant in $M_{\infty}^{\omega}.$
Then we have,
$$\cdots \subset P_{-2} \subset P_{-1} \subset P_0 \subset P_1
 \subset P_2 \subset P_3 \subset \cdots \subset 
P_{\infty} \subset M_{\infty}^{\omega}$$
$$\cdots \subset Q_{-3} \subset Q_{-2} \subset P_1^c \subset P_0^c \subset Q_1
\subset Q_2 \subset \cdots \subset Q_{\infty} \subset M_{\infty}^{\omega}.$$
Actually, we could choose them so that they satisfy
$P_{-k} = Q_k^c$ and $Q_{-k} = P_k^c.$

In the first half of this paper, we aim to show the double commutant 
properties, such as
$P_k^{cc} = P_k$ and $Q_k^{cc}=Q_l$ for $k = 0,1,2, \cdots, \infty.$
When $k < \infty,$ the conclusions follow just from estimating
the Jones indices $[P_k^{cc} : P_k]$ and $[Q_k^{cc} : Q_k].$
It is not so complicated.
However, in case of $k= \infty,$ we could only prove $Q_{\infty}^{cc}
= Q_{\infty}.$  The other one $P_{\infty}^{cc}
= P_{\infty}$ is still open.
The proof deeply owes to the special form of $P_1 = M^{\omega} \cap M'.$
More precisely, when we have a sequence in $P_1$, we actually have
a sequence of sequences of operators.  Then we would like to construct
a new sequence of operators within $P_1$ using these sequences.
Such a construction works for $P_1$, but not for $Q_1$.
%
%$\{F_k\}_{k \in {\bf N}} \subset \omega,$
%for any elements $x^{(k)} = \{ x_n^{(k)} \}_n \in M^{\omega} \cap M'$
%with $k=1,2, \cdots,$
%we construct a new sequence $y := \{ x_n^{(k)} \}_n$ with 
%$n \in F_k \setminus F_{k+1}.$
%Since we may choose each $x_n^{(k)} \in M,$
%it is possible for us to make $y$ in $M^{\omega}.$

In the second half of this paper, 
we eliminate the finite depth condition of $N
\subset M,$ and we shall prove the equivalence of the following conditions.

(1) The subfactor is strongly amenable in Popa's sense.

(2) The central sequence factor $M^{\omega} \cap M'$ has the double
commutant property in $M_{\infty}^{\omega}.$

In the rest, we explain the outline of this proof.
Since $M$ is hyperfinite, we could identify it with
$\otimes_{n=1}^{\infty}M_2({\bf C}),$
and set $A_k : =\otimes_{n=1}^{k}M_2({\bf C}). $
The direction (1) $\Rightarrow$ (2) is easy from 
Ocneanu's central freedom lemma.
For the converse direction (2) $\Rightarrow$ (1), 
we set $A := (M' \cap M_{\infty})' \cap M_{\infty}$
and consider the following non-degenerate commuting squares.
$$\begin{array}{cc}
\begin{array}{ccc}
M & \subset & A \\
\cup &       & \cup \\
A_k' \cap M & \subset & A_k' \cap A 
\end{array} & {\rm with}\, \, k=1,2,3, \cdots
\end{array}
$$
What we want to prove is that the next commuting square is non-degenerate.
$$
\begin{array}{ccc}
M^{\omega} & \subset & A^{\omega} \\
\cup &       & \cup \\
M' \cap M^{\omega} & \subset & M' \cap A^{\omega}. 
\end{array}
$$
Our idea is to ``pile up'' the first commuting squares
and show the non-degeneracy of the second one.
This proof owes much to the paper \cite{PP} as we see in Section 
\ref{kouhan}.

%\end{document}
A part of this work was done at Universit\`a di Roma
``Tor Vergata''.  The author thanks R. Longo for the hospitality.
She is also grateful to 
R. V. Kadison for instructive advice on English,
Y. Kawahigashi for discussions and 
encouragements
and S. Popa for his comments on the first draft of this paper.

\section{Preliminaries}\label{preliminaries}
%A follows from B.$B$H$$$&$N$O(BB$B$G$"$k$3$H$+$i(BA$B$G$"$k$3$H$,$o$+$k!$(B
% denote $BJ*(B by $B5-9f(B $B$@$h!*(B

In the next section we shall study properties of the Jones tower
 arising from a central sequence subfactor
$M' \cap N^ \omega \subset M' \cap M^{\omega}$.  
We use the following notations.
Let $M$ be a hyperfinite II$_1$ factor and $N \subset M$ be 
a type II$_1$ subfactor with finite index and finite depth.
(Here we do not assume the trivial relative commutant condition.)
In the following, we denote the Jones index by $[ M : N ],$ i.e.,
$[ M : N ] = {\rm dim}_N\;L^2(M).$ 
(This definition makes sense whenever $N$ is a II$_1$ factor.
In fact, in Section \ref{kouhan} we define the Jones index without assuming
the factoriality of $M.$)
Then we have the Jones tower 
$$N \subset M \subset M_1 \subset
M_2 \subset M_3 \subset \cdots \subset M_{\infty}.$$
Here $M_{\infty}$ means the weak closure of 
$\bigcup_{k=0}^{\infty}M_{k}$
in the GNS-representation with respect to the trace.\par

%$B$&$k$H$i$W$m$@$/$H(B
In the following, we recall some results about the central sequence subfactor.

(1) 
If we fix a free ultrafilter ${\omega}$ over ${\bf N},$
we obtain an inclusion of type II$_1$ factors;
$M' \cap N^{\omega} \subset M' \cap M^{\omega}$.
We write $M^{\omega}$ for the ultraproduct algebra and
$M' \cap M^{\omega}$ for the central sequence algebra.
The subfactor is called the {\it central sequence subfactor} 
of $N \subset M.$ 
It has been introduced by Ocneanu
\cite{O1}, \cite {O2}.  (See \cite[Section 15.4]{EK}.)\par

(2)
It is known that $M' \cap N^{\omega} \subset M' \cap M^{\omega}$
has a trivial relative commutant.  (See \cite[Lemma 15.25]{EK}.
From the proof of this Lemma, one sees easily that
the trivial relative commutant condition
of $N \subset M$ is not needed for proving that of 
$M' \cap N^{\omega} \subset M' \cap M^{\omega}$!)\par

%$B$0$m!<$P$k$$$s$G$C$/$9(B
(3) 
Let $\{X_i\}_i$ be the set of (isomorphism classes of)
the $N$-$N$ bimodules arising from ${}_N M_M$.  We set
$\gamma=\sum_i [X_i]$ and call it the {\it global index} of
$N\subset M$, where we denote by $[X_i]$ the Jones index of the
bimodule $X_i.$ 
As in \cite[Theorem 12.24, Lemma 15.25]{EK}, we know that 
the index of the asymptotic inclusion of 
$N \subset M$, i.e.,
$M \vee ( M' \cap M_{\infty}) \subset M_{\infty}$, 
is equal to the global index of $N \subset M$, and 
the index of the central sequence 
subfactor is also equal to the global index.  That is,  
$$[ M' \cap M^{\omega} : M' \cap N^{\omega}] = [ M_{\infty} :
 M \vee ( M' \cap M_{\infty})] = \gamma.$$
(This was first noted by Ocneanu \cite{O1}.)
The asymptotic inclusion has been recently studied from many points of view,
such as topological quantum field theory, group actions on factors, 
sector theory, and so on.
When a subfactor is given by a crossed product by an outer action
of a finite group $G$ on the hyperfinite II$_1$ factor $R$, i.e.,
$R \subset R \rtimes G,$
the global index equals to $\# G.$  And
the asymptotic inclusion of $R \subset R \rtimes G$
is given by $Q^{G \times G} \subset Q^G,$
where $Q$
% := R \otimes M_{ \vert G \vert}({\bf C})
%\otimes M_{ \vert G \vert}({\bf C})
%\otimes M_{ \vert G \vert}({\bf C}) \otimes \cdots,
is another hyperfinite II$_1$ factor,
and $G$ is embedded into $G \times G$ by $g \mapsto (g,g).$
Recalling that a paragroup is a generalization of a group,
we may consider that the global index is the ``order'' of a paragroup  
and the asymptotic inclusion gives the quantum double of a paragroup
in an appropriate sense.
(See \cite[Section 12.8]{EK}.)

(4)
The principal graph of the 
asymptotic inclusion of $N \subset M$ is the
connected component of the fusion graph of the bimodule system
arising from $N \subset M$ containing the vertex corresponding to
the bimodule ${}_M M_M.$  Since $N \subset M$ is of finite depth, 
so is $M \vee (M' \cap M_{\infty}) \subset M_{\infty}.$
(See \cite{O2}, \cite[page 663]{EK}.)

(5)
The algebra  $M_{\infty}^{\omega}$ is a II$_1$ factor.  
For simplicity of notations, we set $P_0 := M' \cap N^{\omega}$ and
$P_1 := M' \cap M^{\omega}.$
We denote $ P_{k}' \cap M_{\infty}^{\omega}$ by $P_k^c.$
Then, both $P_0^c$ and $P_1^c$ are II$_1$ factors and
they satisfy $[P_0^c : P_1^c] = [P_1 : P_0] = \gamma.$
(See \cite[Lemmas 15.26, 15.27 and 15.30]{EK}.)

(6)
We also recall the important facts.  From 
the string algebra theory, the central freedom lemma of Ocneanu
(we later explain it in Section \ref{kouhan}) and Popa's generating property, 
we have
$P_0^{cc} = P_0$ and $P_1^{cc} = P_1$.
(This is also due to Ocneanu.  See \cite[Theorem 15.32]{EK}.)
In the next section we study this topic in more detail. 

(7) 
If we choose a downward Jones projection $e \in P_0^c$
with $E_{P_1^c}(e) = 1/ \gamma$ and
set $P_2 := \langle P_1, e \rangle,$
$P_2^c = P_1^c \cap \{ e \} ',$
then we have both $P_{0} \subset P_{1} \subset P_{2}$
and $P_{2}^{c} \subset P_{1}^{c}
 \subset P_{0}^{c}$ are standard.  By repeating this procedure,
we can construct the Jones tower 
$$P_{0} \subset P_{1} \subset P_{2} \subset P_{3} 
\subset \cdots \subset M_{\infty}^{\omega} $$
such that
$$\cdots \subset P_{3}^{c} \subset P_{2}^{c} \subset P_{1}^{c}
 \subset P_{0}^{c}$$
is a tunnel.
 (This has been noted by Ocneanu.  
See \cite[Lemma 15.30]{EK}.)

(8)
Since $P_{1}^{c} \subset P_{0}^{c}$ are type II$_1$ factors 
and $P_0^{cc} = P_0$ and $P_1^{cc} = P_1$ as we have mentioned in (5)
and (6), we also have the Jones tower
$$P_{1}^{c} \subset P_{0}^{c} \subset Q_{1} \subset Q_{2} \subset 
Q_{3} \subset \cdots \subset M_{\infty}^{\omega} $$
such that
$$\cdots \subset Q_{3}^{c} \subset Q_{2}^{c} \subset Q_{1}^{c}
 \subset P_{0} = P_0^{c c} \subset P_{1} = P_1^{c c}$$
is a tunnel. (See \cite[Lemma 15.30]{EK}.) \par

%$B$3$&$5$-$;$s$;$$$N$;$D$a$$(B

We recall the Kosaki index and the Pimsner--Popa index.

Let $A \subset B$ be von Neumann algebras, and $P ( B , A )$ be the
set of all faithful normal semifinite operator-valued weights
from $B$ onto $A$.  In \cite[Theorem 5.9]{H}, 
Haagerup has proved the equivalence between
$P ( B , A ) \neq \emptyset$ and $ P ( A' , B') \neq \emptyset.$
Later, Kosaki
has noticed the existence of
the canonical order-reversing
bijection from $P ( B, A )$ onto $P ( A', B' )$ in \cite{Ko},
which we denote by $E \mapsto E^{-1}$.
When $A, B$ are factors with a fixed normal conditional
expectation $E$, Kosaki defined in \cite[Definition 2.1]{Ko}
${\rm Index}\; E$ as the scalar $E^{-1}(1)$, 
and proved that ${\rm Index}\; E$ does not depend on
a Hilbert space $H$ on which $B$ is represented, i.e., $B \subset B( H )$.
(See \cite[Theorem 2.2]{Ko}.)
When ${\rm Index}\; E  <  \infty $, by setting 
$\tau = (E^{-1}(1))^{-1}$, he also noted that $\tau E^{-1}$ is a 
conditional expectation from $A'$ onto $B'$.\par

Especially when $A$ and $B$ are type II$_1$ factors and
$E : B \to A$ satisfies $\tr _A \circ E = \tr _B$, the Jones index
$[ B : A ] := {\rm dim_A}\;L^2(B) $ equals to ${\rm Index}\; E$.
(See \cite[page 133]{Ko}.)
In the following, we write
$[ B : A ]_{K, E}$ for ${\rm Index}\; E$ to distinguish this index
from the Jones index and the Pimsner--Popa index 
which we define below.

%$B$]$Q(B $B$N$3$H(B
\begin{df}{\bf \cite[Section 2]{PP} \cite[Definition 1.1.1]{P}}
Let $A \stackrel{E }{\subset} B$
 be an inclusion of von Neumann algebras with a conditional expectation
$E $ from $B$ onto $A$. Then we denote
$$(\sup \{ \lambda  \, | \, E (x) \geq \lambda x, x \in B_+ \})^{-1}$$
by $[ B_1 : B_2  ]_{PP, E },$ with the convention $0^{-1} = \infty.$  
\end{df}
In case both $B_1$ and $B_2$ are type II$_1$ factors, this index coincides 
with the Jones index. (See \cite[Proposition 2.1]{PP}.)

\begin{df}{\bf \cite[section 1.1]{P}}
Let $A \stackrel{E}{\subset} B$ be as above.
If $C \stackrel{F}{\subset} D$ is another inclusion
of von Neumann algebras with a conditional expectation
$F$ from $D$ onto $C$  such that $C \subset A, D \subset B$, and
 $F = E|_D$, then we call the square
$$ 
\begin{array}{ccc}
   A & \stackrel{E}{\subset} & B \\
   \cup &                    & \cup \\
   C & \stackrel{F}{\subset} & D
\end{array}$$
a {\rm commuting square}.  If $\overline{{\rm sp}}AD = B,$ we call it
a {\rm non-degenerate} commuting square.
\end{df}
This is a generalization of a case having a trace, which has been
studied in \cite{P0}, \cite{GHJ}.
When the square 
$$\begin{array}{ccc}
   A & \stackrel{E}{\subset} & B \\
   \cup &                    & \cup \\
   C & \stackrel{F}{\subset} & D
\end{array}$$  is a commuting square, we trivially have 
$[ B : A ]_{PP, E} \geq [ D : C ]_{PP, F}$ 
by the definition. \par

The relation between the Kosaki index and the Pimsner--Popa index
is noted in \cite[page 224]{BDH}
and \cite[Theorem 4.1, Corollary 4.2]{L} as follows.

\begin{prop}\label{longo} 
Let $A \stackrel{E}{\subset} B$ be an inclusion of
infinite dimensional factors with a conditional expectation 
$E$ from $B$ onto $A$.
Then we have
$[ B : A ]_{K, E} = [ B : A ]_{PP, E}$.
\end{prop}

\section{Double commutant property }\label{sec1}

In this section we shall prove the double commutant property of 
the Jones tower of $N' \cap M^{\omega} \subset M' \cap M^{\omega}.$
We need the following two easy lemmas. 
\begin{lm} \label{kosaki}
Let $A \stackrel{E}{\subset} B$ 
be  type II$_1$ subfactors of $M_{\infty}^{\omega}$ and we represent these 
von Neumann algebras on $L^2(M_{\infty}^{\omega}).$
Let $E$ be the unique trace-preserving conditional expectation 
from $B$ onto $A$ with ${\rm Index}\; E < \infty$.
 Then we have 
$$[ A' : B' ]_{K, \tau E^{-1}} = [ B : A ].$$
\end{lm}

\begin{proof}
%Since  $B$ is a II$_1$ factor, it has a unique trace $tr$.
%Therefore we have the unique conditional expectation $E : B \to A$
%satisfying $\tr (xy) = \tr ( E(x)y )$ for all $x \in B, y \in A.$ \par 

By the definition of ${\rm Index}\; E$ and the property of the map
$E\mapsto E^{-1}$, we have the following.
\begin{eqnarray*}
[ A' : B' ]_{K, \tau E^{-1}}  &=& {\rm Index}\; (\tau E^{-1})
            = ( \tau E^{-1})^{-1}(1) \\
          &=&  \tau^{-1} (E^{-1})^{-1}(1) = \tau^{-1} 
            = [ B : A ]_{K, E} = [ B : A ]
\end{eqnarray*}

The fourth equality follows from the fact
$(E^{-1})^{-1} = E$.  (See \cite[page 126]{Ko}.) %$B$A$c$s$H$7$i$Y$^$7$g$&(B
\end{proof}

\begin{lm}\label{lemma2}
Let $A \subset B$ and $B^c \subset A^c$
be inclusions of II$_1$ factors contained in
$M_{\infty}^{\omega},$ and represent them on
$L^2(M_{\infty}^{\omega}).$  Let $E$ be the unique trace-preserving
conditional expectation from $A^c$ onto $B^c.$  Then the square
$$\begin{array}{ccc}
A^{c \,\prime} & \stackrel{\tau E^{-1}}{\subset} & B^{c \, \prime} \\
 \cup &                    & \cup \\
A^{c c} & \stackrel{F}{\subset}                 & B^{c c}
\end{array}$$
where $F = (\tau E^{-1})|_{B^{c c}}$
and $\tau = (E^{-1}(1))^{-1},$ is a commuting square.
Furthermore, when we assume that $F(B) \subset A$ and
$A' \cap B = {\bf C},$
$F|_B$ is the unique trace-preserving conditional 
expectation from $B$ onto $A.$ 
\end{lm}
\begin{proof}
Since $M_{\infty}' \subset A^{c \,\prime}\subset B^{c \, \prime}$
and $M_{\infty}'' = M_{\infty},$ the commuting square condition easily follows.
To make it sure, let $x$ and $y$ be arbitrary elements of 
$ B^{c c}$ and $M_{\infty}^{\omega \,\prime}$ respectively.
Since $M_{\infty}^{\omega \,\prime} \subset A^{c \,\prime}
        \subset B^{c \,\prime}, $
$y$ is in $A^{c \,\prime}.$ Then we obtain
$$\tau E^{-1}(x)y = \tau E^{-1}(xy)
= \tau E^{-1}(yx) = y \tau E^{-1}(x),$$
i.e., $$\tau E^{-1}(x) 
\in (M_{\infty}^{\omega \,\prime})' \cap A^{c \,\prime} 
= M_{\infty}^{\omega} \cap A^{c \,\prime} = A^{c c},$$ 
for any $x \in  B^{c c}.$  Therefore $F$ is a conditional
 expectation from $B^{c c}$ onto $A^{c c}$.\par

When $A' \cap B = {\bf C},$ it is known that $A$ has 
the unique conditional expectation onto $B.$  
(See \cite[Proposition 10.17]{S}.)
Thus in this case
$F|_B$ is the unique trace-preserving conditional expectation
from $B$ onto $A.$
%Next, we prove that ${\tr}_A (F(x)) = {\tr}_B (x)$ for any $x \in B.$
%Here, ${\tr}_A$ and ${\tr}_B$ are the unique traces on $A$
%and $B$ respectively.
%By the uniqueness of the trace on $B$, we only have to check if 
%${\tr}_A \circ F$ satisfies the condition of the trace on $B$ or not.
%The positivity and the normality follow trivially.
%The tracial condition has been noted in \cite[page 132]{Ko}.
%And the faithfulness follows from the faithfulness of ${\tr}_A$ and
%$E^{-1} \in P(B^{c \,\prime}, A^{c \,\prime}).$ 
\end{proof}

Now we have two main theorems about the double commutant property.
\begin{th}\label{P_k^{cc}}
We have $P_k^{cc} = P_k$ and $Q_k^{cc} = Q_k$ for any $k \in {\bf N}.$
\end{th}

%by induction$B$H$$$&$N$O7h$^$jJ86g$J$N$GITDj4';l$O$$$i$J$$$s$G$7$?!%(B
\begin{proof}
We prove the theorem by induction on $k$.
We have already known that 
$P_0^{cc} = P_0$ and $P_1^{cc} = P_1$.
Suppose we have $P_k^{cc} = P_k$.
We denote the unique trace-preserving conditional expectation 
from $P_k^{c}$ to $P_{k+1}^{c}$ by $E.$  And we set $\tau$ as above.
We remark that the square
$$ 
\begin{array}{ccc}
   P_k^{c\,\prime} & \stackrel{\tau E^{-1}}{\subset} 
& P_{k+1}^{c\,\prime}\\
   \cup &                    & \cup \\
   P_k^{c c} & \stackrel{F}{\subset} &       P_{k+1}^{c c}
\end{array}$$ where $F = (\tau E^{-1})|_{P_{k+1}^{c c}},$ 
is a commuting square by Lemma \ref{lemma2}.
As we have mentioned at the beginning of this section,
$P_k \subset P_{k+1}$ has a trivial relative commutant.
Thus Lemma \ref{lemma2} also implies that 
$F|_{P_{k+1}}$ is a trace-preserving
conditional expectation from $P_{k+1}$ onto $P_k.$
Applying Lemma \ref{kosaki} and Proposition \ref{longo} to
$P_{k+1}^c \subset P_k^c \subset M_{\infty}^{\omega}$
and using the Pimsner--Popa inequality, we have
\begin{eqnarray*}
\gamma &=& [ P_{k}^c : P_{k+1}^c ]
         = [ P_{k}^c : P_{k+1}^c ]_{K,E}
         = [ P_{k+1}^{c '} : P_{k}^{c '} ]_{K, \tau E^{-1}}
         = [ P_{k+1}^{c '} : P_k^{c '} ]_{PP, \tau E^{-1}}\\
      &\geq & [ P_{k+1}^{c c} : P_k^{c c} ]_{PP, F}
        \geq [ P_{k+1} : P_k ]_{PP, F}
         = [ P_{k+1} : P_k ] = \gamma.
\end{eqnarray*}
Thus we have 
$\gamma = [ P_{k+1} : P_k ]_{PP} = [ P_{k+1}^{c c} : P_k ]_{PP}.$

Since $P_{k+1}^c$ is a II$_1$ factor and $P_{k+1} \subset P_{k+1}^{c c}$
, we obtain
\begin{eqnarray*}
P_{k+1}^{c c} \cap P_{k+1}^{c c '} &\subset & P_{k+1}^{c c} \cap P_{k+1}'
   = ((P_{k+1}^c)' \cap M_{\infty}^{\omega}) \cap P_{k+1}' \\
& = & (P_{k+1}^c)' \cap (P_{k+1}' \cap M_{\infty}^{\omega})
   = (P_{k+1}^c)' \cap P_{k+1}^c 
   = {\bf C}.
\end{eqnarray*}
Hence $P_{k+1}^{c c}$ is a II$_1$ factor. 
%Furthermore, it is a II$_1$ factor 
%because $P_2^{c c}$ is a subalgebra of a II$_1$ factor $M_{\infty}^{\omega}$.
Then we have
$[ P_{k+1} : P_k ]_{PP} = [ P_{k+1} : P_k ]$ and 
$[ P_{k+1}^{c c} : P_k ]_{PP} = [ P_{k+1}^{c c} : P_k ].$
Thus $ [ P_{k+1} : P_k ] = [ P_{k+1}^{c c} : P_k ]$, which means
$P_{k+1}^{c c} = P_{k+1}.$ 

In the same way, we have 
$ [ Q_1^{c c} : P_0^c ]_{PP} = [ Q_1 : P_0^c ]_{PP}$, i.e.,
$Q_1^{cc} = Q_1$.
We can also prove $Q_k^{cc} = Q_k$ in the same way.
\end{proof}

Let $P_{\infty}$ and $Q_{\infty}$ be the weak closures of 
$\bigcup_{k=0}^{\infty}P_{k}$ and $\bigcup_{k=0}^{\infty}Q_{k}$
on $L^2(M_{\infty}^{\omega})$ respectively.  Both 
$P_{\infty}$ and $Q_{\infty}$ are II$_1$ factors, because
$P_k, Q_k$ and $M_{\infty}^{\omega}$ are all II$_1$ factors.
We have the following theorem.

\begin{th} \label{magnolia}
We obtain % $P_{\infty}^{c c} = P_{\infty}$ and 
$Q_{\infty}^{c c} = Q_{\infty}.$
\end{th}
At first glance, this statement may seem strange.  As we mentioned in
Section
\ref{preliminaries},
when we construct $Q_k,$ we have an ambiguity of choosing a Jones 
projection at each step.  
This theorem claims that whatever Jones projections we may
choose, the identity $Q_{\infty}^{c c} = Q_{\infty}$ holds.
In Theorem \ref{P_k^{cc}}, we used only estimates of the Jones index
and nothing particular about ultraproducts was used.  But for this theorem, 
it is essential that $P_1$ is an ultraproduct algebra.
%\begin{lm}\label{clover}
%Let $A_0\supset A_1\supset A_2\supset \cdots$ be a decreasing
%sequence of von Neumann subalgebras of $B\subset C^\omega$
%where $C$ is a finite von Neumann algebra.
%Then we obtain $$( \bigcap_{k=0}^{\infty}A_k)' \cap B
%= \bigvee_{k=0}^{\infty}(A_k' \cap B). $$ \end{lm}

\begin{proof}
We have
\begin{eqnarray*}
Q_{\infty}^{cc} &=& ((\bigvee_{l=1}^{\infty} Q_l)' \cap M_{\infty}^{\omega})'
\cap M_{\infty}^{\omega} \\ &=& ( \bigcap_{l=1}^{\infty}(Q_l' \cap 
M_{\infty}^{\omega}))' \cap M_{\infty}^{\omega} = ( \bigcap_{l=1}^{\infty}
Q_l^c)' \cap M_{\infty}^{\omega},
\end{eqnarray*}
and 
$$Q = \bigvee_{l=1}^{\infty} Q_l = \bigvee_{l=1}^{\infty} Q_l^{cc}
= \bigvee_{l=1}^{\infty} (Q_l^{c \, '} \cap M_{\infty}^{\omega}).$$
It is enough for us to show the following equality.
$$(\bigcap_{l=1}^{\infty}Q_l^c)' \cap M_{\infty}^{\omega} = 
\bigvee_{l=1}^{\infty} (Q_l^{c \, '} \cap M_{\infty}^{\omega})$$ 
It is easy to see that the right hand side is a subalgebra of the 
left hand side.  To prove the converse inclusion, we choose
$$x = \{ x_n\}_n \in (\bigcap_{l=1}^{\infty}Q_l^c)' 
\cap M_{\infty}^{\omega}.$$
If $x$ were not in 
$\bigvee_{l=1}^{\infty} (Q_l^{c \, '} \cap M_{\infty}^{\omega}),$
there would exist $\e > 0$ such that for any $l,$
$$\Vert E_{Q_l^{c \, '} \cap M_{\infty}^{\omega}} (x)-x \Vert_2 > \e.$$
Thus, there exist unitaries $y^l = \{ y_n^l \}_n \in Q_l^c$
such that $\Vert xy^l-y^lx \Vert_2 > \e /2.$
Here we may assume $y_n^l \in M$ because $Q_l^c\subset P_1$.
We denote the Jones projection of the subfactor
$Q_{l+2}^c \subset Q_{l+1}^c$ by $f^l = \{f_n^l\}_n\in Q_l^c.$
%Here we may assume $f_n^l \in M_{\infty }$ because $Q_l^c \subset 
%M_{\infty}^{\omega}.$ 
Let $\{ a_l\}_{l \in {\bf N}}$ be an $L^2$-dense subset of $M.$  
Then we have
$$\begin{array}{ccccl}
f^3&\in& Q_3^c&=&Q_2^c\cap\{f^1\}'\\
&&\cap&&\\
f^2&\in& Q_2^c&=&Q_1^c\cap\{f^0\}'\\
&&\cap&&\\
f^1&\in& Q_1^c&=&P_0\cap\{f^{-1}\}'\\
&&\cap&&\\
f^0&\in& P_0&=&P_1\cap\{f^{-2}\}'\\
&&\cap&&\\
f^{-1}&\in& P_1&=&M'\cap M^\omega\\
&&\cap&&\\
f^{-2}&\in& P_2&=&(M'\cap M^\omega)\vee\{f^{-2}\}\\
&&\cap&&\\
&&\vdots&&\\
&&\cap&&\\
&&M_{\infty}^{\omega}&&
\end{array}$$ 
which means
$$y^l \in Q_l^c = M^{\omega} \cap M' \cap \{f^{-2}, f^{-1}, 
f^0, \cdots , f^{l-2} \} '.$$

Let $F_0 := {\bf N}$ and
$$ F_k := F_{k-1} \cap [k, \, \infty ) \cap
\left\{
\begin{array}{c|l}
&  \Vert x_ny_n^k-y_n^kx_n \Vert_2 > \e /2 \\
n &  \Vert f_n^iy_n^k-y_n^kf_n^i \Vert_2 < 1/k,\ 
 {\rm for}\  i = -2, -1, \cdots, k-2 \\
&   \Vert a_iy_n^k-y_n^ka_i \Vert_2 < 1/k, \
{\rm for}\  i = 1, 2, \cdots, k
\end{array}
\right\}. $$
Then, each $F_k$ is in $\omega $ and $\bigcap_{k=1}^{\infty} F_k = \emptyset.$
If we set $y := \{ y_n^k \}_n$ for $n \in F_k \setminus F_{k+1},$ then
$$y \in M^{\omega} \cap M' \cap \{ f^{-2}, f^{-1}, f^0, \cdots \} '
= \bigcap_{k=1}^{\infty} Q_k^c,$$
and $\Vert xy-yx \Vert_2 > \e /2.$  Thus, 
$x \notin ( \bigcap_{k=1}^{\infty} Q_k^c)' \cap M_{\infty}^{\omega},$
which is a contradiction.
\end{proof}

Unfortunately, we cannot prove $P_{\infty}^{cc} = P_{\infty}$
in the same way.  In the above proof,
we choose every $y_n^k$ in $M,$ thus $y = \{y_n^k \}_n \in M^{\omega}.$
However, since $P_0^c$ is represented as
$P_0^c = \bigvee_{k=1}^{\infty} (A_{0 \, \infty} \vee
A_{0 \, k})^{\omega},$ (here the string algebra $A_{k,l}$ does not matter,
see \cite[Section 15]{EK} for more details),
if we construct an element of the filter and $y$ as above, we are not sure
whether such $y$ could be in $P_0^c$ or not.

%$BDjM}$r=q$/$H$-$K(B
%\begin{theorem}\label{abc}...\end{theorem}
%$B$H=q$1$P!$$3$NDjM}$NL>A0$O(Babc$B$K$J$j$^$9!%(B($B%"%k%U%!%Y%C%H$G;O$^$l$P$I$s$J(B
%$BL>A0$G$b$D$1$i$l$^$9!%(B)$B$=$l$r0zMQ$9$k$H$-$O!$(BTheorem \ref{abc}$B$N$h$&$K$7$^$9!%(B

\section{Applications to paragroups}\label{application}

In this section we study the double sequences of 
the higher relative commutants of the subfactor
$M' \cap N^{\omega} \subset M' \cap
M^{\omega}.$
Owing to the double commutant properties, we can see the relations
among the higher relative commutants clearly.

\begin{lm}\label{cherry}
We have the following identities.

(1) $\bigvee_{k=0}^{\infty}(P_k \cap Q_l) = P_{\infty} \cap Q_l$

(2) $\bigvee_{l=0}^{\infty}(P_k \cap Q_l) = P_k \cap Q_{\infty}$

(3) $\bigvee_{k=0}^{\infty}(P_k \cap Q_{\infty}) = P_{\infty} \cap Q_{\infty}$

(4) $\bigvee_{l=0}^{\infty}(P_{\infty} \cap Q_l) = P_{\infty} \cap Q_{\infty}$
\end{lm}

\begin{proof}
(1) It is clear that the left hand side is contained in the right hand side.
We prove the converse inclusion.
%Let $x$ be an element of $P_{\infty} \cap Q_l.$
%We set $x_k := E_{P_k \cap Q_l}(x).$
Since $Q_l^{cc} = Q_l$ and $Q_l^c \subset P_k \subset P_{\infty},$ the square
$$\begin{array}{ccc}
P_k \cap Q_l &\subset & P_{\infty} \cap Q_l \\
   \cap &                    & \cap \\
P_k     & \stackrel{E}{\subset} & P_{\infty}
\end{array}$$  is a commuting square.
Then for any $x \in P_{\infty} \cap Q_l,$ we have
$E_{P_k \cap Q_l}(x) = E_{P_k}(x)$ and 
$\Vert E_{P_k \cap Q_l}(x)  \Vert_{\infty} 
\leq \Vert x \Vert_{\infty}.$
Since $\bigvee_{k=1}^{\infty}P_k = P_{\infty},$
we have $$\Vert x-E_{P_k \cap Q_l}(x) \Vert_2 
= \Vert x-E_{P_k}(x) \Vert_2 \to 0,$$
which means $E_{P_k \cap Q_l}(x)$ converges to $x$ strongly.
Therefore, we have
$E_{\vee_{k=1}^{\infty} (P_k \cap Q_{\infty})}(x) = x$ for any 
$x \in P_{\infty} \cap Q_{\infty}.$

(2), (3) Since $P_k^{cc} = P_k$ and $Q_{\infty }^{cc} = Q_{\infty},$ 
then the squares
$$\begin{array}{ccc} 
P_k \cap Q_l & \subset & Q_l \\
\cap & & \cap \\
P_k \cap Q_{\infty }& \subset & Q_{\infty}
\end{array} $$ and
$$\begin{array}{ccc} 
P_k \cap Q_{\infty} & \subset & P_{\infty} \cap Q_{\infty} \\
\cap & & \cap \\
P_k & \subset & P_{\infty} 
\end{array}$$ are commuting squares.  Thus we obtain equalities
(2), (3) in the same way as is the proof of (1).

(4) The equalities below show (4).
\begin{eqnarray*}
\bigvee_{l=1}^{\infty} (Q_l \cap P_k) &=& 
\bigvee_{l=1}^{\infty} (Q_l^{c \, '} \cap P_k)
= ( \bigcap_{l=1}^{\infty} Q_l^c)' \cap P_k 
= ( \bigcap_{l=1}^{\infty} Q_l^c)' \cap M_{\infty}^{\omega} \cap P_k \\
&=& (\bigvee_{l=1}^{\infty} (Q_l^{c \, '} \cap M_{\infty}^{\omega})) \cap P_k
= (\bigvee_{l=1}^{\infty} Q_l^{cc}) \cap P_k = Q_{\infty} \cap P_{\infty}.
\end{eqnarray*}
Both the second and fourth identities follow from the same arguments as
in the proof of Theorem \ref{magnolia}.
\end{proof}

In addition to the properties mentioned in Section \ref{preliminaries},
the subfactor $N^{\omega}\cap M'\subset M' \cap M^{\omega}$
is known to have the finite depth property.  To prove it, the next
theorem of Ocneanu has been very useful.

\begin{th}{\bf (Ocneanu)}\label{dual}
Let $N \subset M$ be a subfactor of the hyperfinite II$_1$ factor
with finite index and finite depth.  The paragroups of the 
central sequence subfactor $N^{\omega}\cap M'\subset M' \cap M^{\omega}$
and of the asymptotic inclusion 
$M\vee (M'\cap M_{\infty})\subset M_{\infty}$
are mutually dual.
\end{th}

This has been noted by Ocneanu in \cite{O2}.  One can see a proof
in \cite[Theorem 15.32]{EK}.
In general, we say that two paragroups are {\it dual} to each other
if and only if the corresponding subfactors are dual to each other.
(See \cite[page 570]{EK}.)
Originally the adjective  ``dual'' comes from the following duality of 
groups.
When $G$ is a finite group and $R$ a II$_1$ factor, 
we have a fixed point algebra and a crossed product
by an outer action, i.e., $R^G \subset R$ and
$R \subset R \rtimes G.$  It is known that the
$R^G$-$R^G$ bimodules are indexed by $\hat G$ and the
$R$-$R$ bimodules are indexed by $G.$ 
Thus we say $R^G \subset R$ and $R \subset R \rtimes G$
are dual.  Extending this duality, we say
$N \subset M$ and $M \subset M_1$ are dual
when $N \subset M \subset M_1$ is standard.

The finite depth condition of $N^{\omega}\cap M'\subset M' \cap M^{\omega}$
follows from this theorem and the finite depth condition of 
$M\vee (M'\cap M_{\infty})\subset M_{\infty}.$
Thanks to the above arguments,
especially using the double commutant properties
($P_k^{cc} = P_k$ and $Q_l^{cc} = Q_l$), we simplify the proof of
\cite[Theorem 15.32]{EK}.  That is,
in \cite{EK} they have shown $P_0' \cap P_k = (P_k^c)' \cap P_0^c$
by using several inclusions and two anti-isomorphisms.
However, with the double commutant properties, it is quite natural for us
to write the higher relative commutants by the combination of
algebras, one from 
$$\cdots \subset Q_{3}^{c} \subset Q_{2}^{c} \subset Q_{1}^{c}
 \subset P_{0} \subset P_{1} \subset P_{2} \subset P_{3} 
\subset \cdots \subset M_{\infty}^{\omega} $$
and the other from
$$\cdots \subset P_{3}^{c} \subset P_{2}^{c} \subset P_{1}^{c}
 \subset P_{0}^{c} \subset Q_{1} \subset Q_{2} \subset 
Q_{3} \subset \cdots \subset M_{\infty}^{\omega}.$$

\begin{proof}
The double sequence of the higher relative commutants of the 
central sequence subfactor $N^{\omega}\cap M'\subset M' \cap M^{\omega}$ is 
given as in the following diagram.  Here we note that 
$N^{\omega}\cap M'\subset M' \cap M^{\omega}$ has a trivial relative commutant,
(see \cite[Lemma 15.25]{EK}), and use the conventions of Lemma \ref{cherry}.
$$\begin{array}{ccccccccccccccc}
&&&&&&{\bf C} & \!\!\!\subset\!\!\! &  {\bf C} & \!\!\!\subset\!\!\! &  P_3\cap P_1^c
& \!\!\!\subset\!\!\! & \cdots & \!\!\!\subset\!\!\! &  P_{\infty} \cap P_1^c \\
&&&&&&\cap && \cap && \cap &&&& \cap \\
&&&&{\bf C}&\!\!\!\subset\!\!\! &  {\bf C}& \!\!\!\subset\!\!\! &  P_2\cap P_0^c 
& \!\!\!\subset\!\!\! &  P_3\cap P_0^c & \!\!\!\subset\!\!\! & \cdots & \!\!\!\subset\!\!\! & 
P_{\infty}\cap P_0^c \\
&&&&\cap && \cap && \cap && \cap &&&& \cap \\
&& {\bf C}&\!\!\!\subset\!\!\! &  {\bf C}&\!\!\!\subset\!\!\! &  P_1\cap Q_1\ 
&\!\!\!\subset\!\!\! &  P_2\cap Q_1 & \!\!\!\subset\!\!\! & 
P_3\cap Q_1 &\!\!\!\subset\!\!\! & \cdots & \!\!\!\subset\!\!\! & P_{\infty}\cap Q_1 \\
&&\cap && \cap && \cap && \cap && \cap& &&  & \cap \\
{\bf C} &\!\!\!\subset\!\!\! &  {\bf C}& \!\!\!\subset\!\!\! &P_0\cap Q_2 &\!\!\!\subset\!\!\! &P_1\cap Q_2 
& \!\!\!\subset\!\!\! & P_2\cap Q_2 & \!\!\!\subset\!\!\! &  P_3 \cap Q_2 & \!\!\!\subset\!\!\! & \cdots & \!\!\!\subset\!\!\! 
& P_{\infty}\cap Q_2 \\
\cap & &\cap & &\cap & &\cap & &\cap & &\cap &&&&\cap\\
\vdots &&\vdots &&\vdots &&\vdots &&\vdots &&\vdots &&&&
\vdots \\
\cap & &\cap & &\cap & &\cap & &\cap & &\cap &&&&\cap\\
Q_2^c\cap Q_{\infty}&\!\!\!\subset\!\!\! &Q_1^c\cap Q_{\infty}&\!\!\!\subset\!\!\! &
P_0\cap Q_{\infty}&\!\!\!\subset\!\!\! &P_1\cap Q_{\infty}&\!\!\!\subset\!\!\! &
P_2\cap Q_{\infty}&\!\!\!\subset\!\!\! &P_3\cap Q_{\infty}&\!\!\!\subset\!\!\! &\cdots
&\!\!\!\subset\!\!\! &P_{\infty}\cap Q_{\infty}\\
\end{array}$$

On the other hand,
the double sequence of the higher relative commutants of 
$M\vee (M'\cap M_{\infty})\subset M_{\infty}$ is equal to
that of $P_1^c \subset P_0^c.$  (See \cite[Lemma 15.31]{EK}.)  
Then the double sequence is as follows.

$$\begin{array}{ccccccccccccccc}
&&&&&&{\bf C} & \!\!\!\subset\!\!\! &  {\bf C} & \!\!\!\subset\!\!\! &  Q_2\cap P_0
& \!\!\!\subset\!\!\! & \cdots & \!\!\!\subset\!\!\! &  Q_{\infty} \cap P_0 \\
&&&&&&\cap && \cap && \cap &&&& \cap \\
&&&&{\bf C}&\!\!\!\subset\!\!\! &  {\bf C}& \!\!\!\subset\!\!\! &  Q_1\cap P_1
& \!\!\!\subset\!\!\! &  Q_2\cap P_1 & \!\!\!\subset\!\!\! & \cdots & \!\!\!\subset\!\!\! & 
Q_{\infty}\cap P_1 \\
&&&&\cap && \cap && \cap && \cap &&&& \cap \\
&& {\bf C}&\!\!\!\subset\!\!\! &  {\bf C}&\!\!\!\subset\!\!\! &  P_0^c\cap P_2\ 
&\!\!\!\subset\!\!\! &  Q_1\cap P_2 & \!\!\!\subset\!\!\! & 
Q_2\cap P_2 &\!\!\!\subset\!\!\! & \cdots & \!\!\!\subset\!\!\! & Q_{\infty}\cap P_2 \\
&&\cap && \cap && \cap && \cap && \cap& &&  & \cap \\
{\bf C} &\!\!\!\subset\!\!\! &  {\bf C}& \!\!\!\subset\!\!\! &P_1^c\cap P_3 &\!\!\!\subset\!\!\! &P_0^c\cap P_3 
& \!\!\!\subset\!\!\! & Q_1\cap P_3 & \!\!\!\subset\!\!\! &  Q_2 \cap P_3 & \!\!\!\subset\!\!\! & \cdots & \!\!\!\subset\!\!\! 
& Q_{\infty}\cap P_3 \\
\cap & &\cap & &\cap & &\cap & &\cap & &\cap &&&&\cap\\
\vdots &&\vdots &&\vdots &&\vdots &&\vdots &&\vdots &&&&
\vdots \\
\cap & &\cap & &\cap & &\cap & &\cap & &\cap &&&&\cap\\
P_3^c\cap P_{\infty}&\!\!\!\subset\!\!\! &P_2^c\cap P_{\infty}&\!\!\!\subset\!\!\! &
P_1^c\cap P_{\infty}&\!\!\!\subset\!\!\! &P_0^c\cap P_{\infty}&\!\!\!\subset\!\!\! &
Q_1\cap P_{\infty}&\!\!\!\subset\!\!\!   &Q_2\cap P_{\infty}&\!\!\!\subset\!\!\! &\cdots
&\!\!\!\subset\!\!\! &Q_{\infty}\cap P_{\infty} 
\end{array}
$$ 

Since the asymptotic inclusion is anti-isomorphic to itself, its
paragroup is opposite to itself, so the above double sequence is
isomorphic to the following.

$$\begin{array}{ccccccccccccc}
&&&&&&{\bf C} & \!\!\!\subset\!\!\! &  {\bf C} 
& \!\!\!\subset\!\!\! & \cdots & \!\!\!\subset\!\!\! &  P_{\infty} \cap P_2^c \\
&&&&&& \cap && \cap &&&& \cap \\
&&&&{\bf C} & \!\!\!\subset\!\!\! &  {\bf C} & \!\!\!\subset\!\!\! &  P_3\cap P_1^c
& \!\!\!\subset\!\!\! & \cdots & \!\!\!\subset\!\!\! &  P_{\infty} \cap P_1^c \\
&&&&\cap && \cap && \cap &&&& \cap \\
&&{\bf C}&\!\!\!\subset\!\!\! &  {\bf C}& \!\!\!\subset\!\!\! &  P_2\cap P_0^c 
& \!\!\!\subset\!\!\! &  P_3\cap P_0^c & \!\!\!\subset\!\!\! & \cdots & \!\!\!\subset\!\!\! & 
P_{\infty}\cap P_0^c \\
&&\cap && \cap && \cap && \cap &&&& \cap \\
{\bf C}&\!\!\!\subset\!\!\! &  {\bf C}&\!\!\!\subset\!\!\! &  P_1\cap Q_1\ 
&\!\!\!\subset\!\!\! &  P_2\cap Q_1 & \!\!\!\subset\!\!\! & 
P_3\cap Q_1 &\!\!\!\subset\!\!\! & \cdots & \!\!\!\subset\!\!\! & P_{\infty}\cap Q_1 \\
\cap && \cap && \cap && \cap && \cap& &&  & \cap \\
{\bf C}& \!\!\!\subset\!\!\! &P_0\cap Q_2 &\!\!\!\subset\!\!\! &P_1\cap Q_2 
& \!\!\!\subset\!\!\! & P_2\cap Q_2 & \!\!\!\subset\!\!\! &  P_3 \cap Q_2 & \!\!\!\subset\!\!\! & \cdots & \!\!\!\subset\!\!\! 
& P_{\infty}\cap Q_2 \\
\cap & &\cap & &\cap & &\cap & &\cap &&&&\cap\\
\vdots &&\vdots &&\vdots &&\vdots &&\vdots &&&&
\vdots \\
\cap & &\cap & &\cap & &\cap & &\cap &&&&\cap\\
Q_1^c\cap Q_{\infty}&\!\!\!\subset\!\!\! &
P_0\cap Q_{\infty}&\!\!\!\subset\!\!\! &P_1\cap Q_{\infty}&\!\!\!\subset\!\!\! &
P_2\cap Q_{\infty}&\!\!\!\subset\!\!\! &P_3\cap Q_{\infty}&\!\!\!\subset\!\!\! &\cdots
&\!\!\!\subset\!\!\! &P_{\infty}\cap Q_{\infty} 
\end{array}$$

By shifting the first diagram by one line vertically, we get the third
diagram, which means the paragroup of $N^{\omega} \cap M' \subset M' \cap M^{\omega}$
is dual to that of $M \vee (M' \cap M_{\infty}) \subset M_{\infty}.$
\end{proof}

\begin{th}\label{peppermint}
The von Neumann algebra $P_{\infty} \cap Q_{\infty}$
is a hyperfinite II$_1$ factor.
\end{th}

\begin{proof}
%We first show that $P_k \cap Q_{\infty}$ is a hyperfinite II$_1$ factor.
Since $$P_k^c \subset \cdots \subset P_1^c \subset P_0^c \subset Q_1
\subset Q_2 \subset Q_3 \subset \cdots $$
is the Jones tower with $[ P_0^c : P_1^c ] = \gamma < \infty$, we have 
${\rm dim}(Q_l \cap P_k) = {\rm dim}(Q_l \cap P_k^{c \,\prime}) < \infty.$
Thus Lemma \ref{cherry} and the finite depth property of $P_1^c \subset P_0^c$ 
imply that $P_{\infty} \cap Q_{\infty}$
is a hyperfinite II$_1$ factor.
%$P_k \cap Q_{\infty}$ is a hyperfinite II$_1$ factor.
%Lemma \ref{fin.depth} implies that
%$$Q_1^c \cap Q_{\infty} \subset P_0 \cap Q_{\infty}
%\subset P_1 \cap Q_{\infty} \subset P_2\cap Q_{\infty}
%\subset P_3 \cap Q_{\infty} \cdots$$ 
%is the Jones tower given by the string algebra construction 
%from the paragroup of the subfactor $P_1^c \subset P_0^c.$
%Thus we have proved that $P_{\infty} \cap Q_{\infty}$ is a hyperfinite
%II$_1$ factor by Lemma \ref{cherry}.
\end{proof} 

\begin{cor}
The II$_1$ factors $P_{\infty}$ and $Q_{\infty}$ are strictly 
smaller than $M_{\infty}^{\omega}.$
\end{cor}

By the definition of $P_0,$ it is not separable and quite a large algebra. 
Thus $P_{\infty}$ is extremely large.  
Once we take an ultraproduct of
$M_{\infty},$ it becomes extraordinarily as large as it contains 
$P_{\infty}$ strictly. 

\begin{proof}
Suppose $P_{\infty} = M_{\infty}^{\omega}$ on the contrary,
then we have $P_{\infty} \cap Q_{\infty} = M_{\infty}^{\omega}
 \cap Q_{\infty} = Q_{\infty}.$
Since $Q_{\infty}$ contains $P_0^c$ which is not hyperfinite, 
$Q_{\infty}$ is not hyperfinite.  By Lemma \ref{peppermint},  
$P_{\infty} \cap Q_{\infty}$ is 
hyperfinite, which is a contradiction.  
We remark that if we suppose $Q_{\infty} = M_{\infty}^{\omega},$
we similarly obtain a contradiction.
\end{proof}

\section{The strong amenability and the central sequence factor 
$M^{\omega} \cap M'$}
\label{kouhan}

In this section we shall prove the following theorem.

\begin{th}\label{ame}
Let $M$ be a hyperfinite II$_1$ factor and $N \subset M$ be 
a type II$_1$ subfactor with finite index.
(Here we do {\rm not} assume the trivial relative commutant condition
 nor {\rm finite depth} condition.)
The following are equivalent.\par
(1)$(M' \cap M_{\infty})' \cap M_{\infty} = M.$ \par
(2)$(M_{\omega}' \cap M_{\infty}^{\omega})' \cap M_{\infty}^{\omega} =
 M_{\omega}.$
\end{th}

We remark condition (1) is equivalent to the strong amenability of Popa.
(See \cite[Theorem 4.1.2]{Pacta}.)
He has proved that the strong amenability is equivalent to
his generating property.
When a subfactor has a finite depth, it is easy to see that it is
strongly amenable.
In Section \ref{sec1}, we have assumed the finite depth condition of 
$N \subset M,$ therefore condition (1) always holds.

Before the proof,
we recall Ocneanu's central freedom lemma (\cite{O1}, \cite[Lemma 15.20]{EK})
as follows, since we use it here for several times.

\begin{lm}{\bf (Central freedom lemma, Ocneanu)}
Let $L \subset P \subset Q$ be finite von Neumann algebras
and $L$ a hyperfinite factor.  Then we obtain
$$(L' \cap P^{\omega})' \cap Q^{\omega} = L \vee (P' \cap Q)^{\omega}.$$
\end{lm}

In this lemma, the hyperfiniteness condition is indispensable
because we approximate $L$ by $L_m := \otimes_{i=1}^m M_2({\bf C})$
and use the finite dimensionality of $L_m$ essentially.

\begin{proof} {\bf (Theorem \ref{ame} (1) $\Rightarrow$ (2))}
By using the central freedom lemma twice it is straightforward to see
that (2) follows from (1), i.e.,
\begin{eqnarray*}
(M_{\omega}' \cap M_{\infty}^{\omega})' \cap M_{\infty}^{\omega} &=&
((M' \cap M^{\omega})' \cap M_{\infty}^{\omega})' \cap M_{\infty}^{\omega}
= (M \vee (M' \cap M_{\infty})^{\omega})' \cap M_{\infty}^{\omega}\\
&=& M' \cap (M' \cap M_{\infty})^{\omega \,\prime} \cap M_{\infty}^{\omega}
=  M' \cap ((M' \cap M_{\infty})' \cap M_{\infty})^{\omega}\\
&=& M' \cap M^{\omega}.
\end{eqnarray*} \end{proof}

We prove the converse direction in the following way.
$$
((2) \Rightarrow (1))
\stackrel{(\ref{omeg})}{\Leftarrow} 
(\ref{AM})
\left\{
\begin{array}{l}
\begin{array}{l}
{\rm The\ case\ }\left[A^{\omega} : M^{\omega}\right] < \infty  \\
(\ref{PPbasis})
\end{array}\\
\begin{array}{l}
{\rm The\ case\ }\left[A^{\omega} : M^{\omega}\right] = \infty \\
(\ref{infty}) (\ref{dim})
\end{array}
\left\{
\begin{array}{l}
{\rm The\ case\ }\dim (M' \cap A) = \infty \\
                           \\
{\rm The\ case\ }\dim (M' \cap A) < \infty  \\
(\ref{maximal}) (\ref{kityaku}) \\
\end{array}
\right.\\
\end{array}
\right.
$$

Though we do not know whether $A$ is a factor or not, 
$[A : M]$ has a meaning
in the sense of Section 
\ref{preliminaries}, i.e., $[A : M] := \dim_M L^2(A)$.
We denote %the integer part of $[A : M] \in [1, \infty]$ by $n,$
the unique trace-preserving conditional expectation from
$M_{\infty}$ onto $M$ by $E_M.$  \par

\begin{prop}\label{infty}
Let $A$ be a finite von Neumann algebra with a fixed trace $\tr,$ and
$M$ its type II$_1$ subfactor.  If,
for any $\e >0,$ there exists a non-zero projection $p \in A$ such that
$E_M (p) \leq \e 1_M,$ then we obtain
$[A : M]_{PP} = \infty.$
\end{prop}

In \cite[page 71]{PP}, this proposition 
has been used since it is trivial from
the definition of the Pimsner--Popa index.
We include a proof here since $A$ is not a factor now,
but actually this does not cause any trouble.
\begin{proof}
We suppose $[A : M]_{PP} < \infty$ on the contrary, then
there exists $\e_0 >0$ such that for any $x \in A_+$ we have
$E_M (x) \geq 2\e_0 x.$  From the assumption 
there exists a projection $0 \neq p_0 \in A$
such that $E_M (p) \leq \e_0 1_M.$
Then we have $2 \e_0 p_0 \leq E_M (p) \leq \e_0 1_M.$
If we multiply all the three operators by $p_0$ both from the right and 
left, we obtain
 $2 \e_0 p_0 \leq p_0 E_M (p)p_0 \leq \e_0 p_0,$
which is a contradiction.
\end{proof}

We also need the following easy lemma, so we remark it here.
\begin{prop}\label{dim}
Let $A$ be a finite von Neumann algebra with a fixed trace $\tr,$ and
$M$ its II$_1$ subfactor. Then we have the following.

(1) When $\dim (M' \cap A) < \infty,$
we have the following identity.
$$[A : M]_{PP} = \max_{1 \leq i \leq k} \frac{[Aq_i : Mq_i]}{\tr (q_i)}$$
In particular,
the condition $[A : M]_{PP} = \infty$ is equivalent to
the condition $[A : M] = \infty.$

(2) When $\dim (M' \cap A) = \infty,$ we obtain
$[A : M]_{PP} = \infty = [A: M].$
\end{prop}
\begin{proof}
(1) Since $\dim (A' \cap A) \leq \dim (M' \cap A) < \infty,$
there exist a finite number of minimal central projections 
$q_1, \cdots, q_k \in A,$
such that $A = Aq_1 \oplus \cdots \oplus Aq_k.$
If we set the isomorphism $\phi_i : Mq_i \to M$ by $\phi_i (yq_i) := y,$
for any $x = x_1 + \cdots + x_k \in \oplus_{i=1}^k Aq_i,$ we have
\begin{eqnarray*}
E_M (x) &=& \sum_{i=1}^k E_M (x_iq_i)
= \sum_{i=1}^k \phi_i (E_{Mq_i} (x_i)) \tr(q_i) \\
&=& \sum_{i=1}^k \tr (q_i) \sum_{j=1}^k \phi_i (E_{Mq_i} (x_i)) q_j
\geq \sum_{i=1}^k \tr (q_i) E_{Mq_i} (x_i) \\
&\geq & \sum_{i=1}^k  \frac{ \tr (q_i)x_i}{[Aq_i : Mq_i]} 
\geq  \left( \max_{1 \leq i \leq k} \frac{[Aq_i : Mq_i]}{\tr (q_i)}
\right)^{-1} x.
\end{eqnarray*}
Thus we obtain 
$$[A : M]_{PP} \leq \max_{1 \leq i \leq k} \frac{[Aq_i : Mq_i]}{\tr (q_i)}$$

In case $[A : M]_{PP} = \infty,$
the above inequality is enough for us, because the following identities 
\begin{eqnarray*}
[A : M] &=& \sum_{i=1}^k \dim_M L^2 (q_i A q_i) \\
&=& \sum_{i=1}^k \dim_{Mq_i} L^2 (q_i A q_i) 
= \sum_{i=1}^k [q_i A q_i : Mq_i],  \end{eqnarray*}
imply $[A : M]= \infty.$

In case $[A : M]_{PP} < \infty,$ 
for any $i \in \{ 1, \cdots, k \}$ and $x \in A,$ 
we have the following inequality.
$$E_M(xq_i ) \geq \frac{xq_i}{[A : M]_{PP}}$$
Because the left hand side is equal to
$\tr (q_i) \sum_{j=1}^k \phi_i (E_{Mq_i} (x_i)) q_j,$
we obtain
$$\tr (q_i) E_{Mq_i} (x_i) \geq \frac{xq_i}{[A : M]_{PP}}.$$
Thus  $[Aq_i : Mq_i]_{PP} < \infty.$
Since  $[Aq_i : Mq_i]_{PP}= [Aq_i : Mq_i],$
we have
$[A q_i : Mq_i] < \infty.$
Let $i_0$ be an index with
$$\max_{1 \leq i \leq k} \frac{[Aq_i : Mq_i]}{\tr (q_i)} =
\frac{[Aq_{i_0} : Mq_{i_0}]}{\tr (q_{i_0})}.$$
Let $e_i \in Aq_i$ be a Jones projection of $Mq_i \subset Aq_i$ for
the downward basic construction.
Then we have 
$$E_M (e_{i_0}) = 
\left( \max_{1 \leq i \leq k} \frac{[Aq_i : Mq_i]}{\tr (q_i)} \right)^{-1}
\geq 
\left( \max_{1 \leq i \leq k} \frac{[Aq_i : Mq_i]}
{\tr (q_i)} \right)^{-1} e_{i_0}.$$
Thus we obtain
$$[A : M]_{PP} = \max_{1 \leq i \leq k} \frac{[Aq_i : Mq_i]}{\tr(q_i)}.$$

The identity
$$[A : M] = \sum_{i=1}^k [q_i A q_i : Mq_i]$$
implies that $[A : M]= \infty$ is equivalent to $[A : M]_{PP}= \infty.$

(2) Since $A$ is a finite von Neumann algebra and 
$\dim (M' \cap A) = \infty,$
for any $\e > 0,$ there exists a non-zero projection $e \in M' \cap A$
such that $0 \neq \tr (e) \leq \e.$
Since the square$$
\begin{array}{ccc}
M & \subset & A \\
\cup &     &  \cup \\
M'\cap M = {\bf C} & \subset & M' \cap A
\end{array} $$
is a commuting square, we have
$0 \neq E_M (e) = \tr (e) \leq \e 1_M.$
Lemma \ref{infty} implies the first equality
$[A : M]_{PP} = \infty.$  \par

In the rest of the proof (2), 
the idea is given by \cite[page 71]{PP}.
They use Jones' identity.  In our case $A$ is not a II$_1$ factor,
then we cannot use it in his original form. Thus we replace
the identity  as follows.

We suppose $\dim (M' \cap A) = \infty$ and $[A : M] < \infty.$
For any non-zero projection $p \in M' \cap A,$ we have
\begin{eqnarray*}
[Ap : Mp] &=& \dim_{Mp} L^2(pAp) \\
&=& \dim_M L^2(A) \tr_{M'}(p) \tr_{M'}(p').
\end{eqnarray*}
Here $p'$ means right multiplication of $p.$ We remark that
$M'$ is a II$_1$ factor by the assumption $[A : M] < \infty.$
For any $m,$ there exist projections $p_1, \cdots, p_m \in M' \cap A$
such that $\sum_{i=1}^m p_i =1.$  Then we have
$[A : M] = \sum_{i=1}^m [Ap_i : Mp_i]/ \tr (p_i) \geq m.$
Therefore $[A : M] = \infty,$ which is a contradiction.
\end{proof}

We need the next two lemmas to show the
invariance of the Jones indices under taking ultraproducts.
\begin{lm}\label{omega} {\bf \cite[Proposition 1.10]{PP}}
Let ${\cal N} \subset {\cal M}$ be 
II$_1$ factors. Then $[{\cal M}^{\omega} : {\cal N}^{\omega}] =
 [{\cal M} : {\cal N}]$
\end{lm}

This identity also holds in the case of the infinite index.
Thanks to the above lemma, the next one is quite natural,
where we have dropped the factoriality of $Q.$
\begin{lm}\label{omeg}
Let $Q$ be a finite von Neumann algebra with a fixed trace $\tr,$ and
$P$ be a II$_1$ factor.
We have $[Q : P] = [Q^{\omega} : P^{\omega}].$
%If $[Q^{\omega} : P^{\omega}] = 1,$ then $[Q : P] = 1,$ i.e., $Q=P.$ 
\end{lm}
\begin{proof}
When $\dim (P' \cap Q) = \infty,$
by the central freedom lemma, we have 
$P^{\omega \, \prime} \cap Q^{\omega} = (P' \cap Q)^{\omega}$
then $\dim (P^{\omega \, \prime} \cap Q^{\omega}) = \infty.$ 
Thus $[Q^{\omega} : P^{\omega}] = \infty = [Q : P]$
by Proposition \ref{dim}.

When $\dim (P' \cap Q) < \infty,$
there exist minimal central projections $q_1, \cdots, q_n \in Q$
such that $Q= \oplus_{i=1}^n Qq_i.$
Then we have
\begin{eqnarray*}
[Q^{\omega} : P^{\omega}] &=& [ \oplus_{i=1}^n Q^{\omega}q_i :
\oplus_{i=1}^n P^{\omega}q_i] = \sum_{i=1}^n [Q^{\omega}q_i :
P^{\omega}q_i] \\
&=& \sum_{i=1}^n [Qq_i : Pq_i] = [Q : P]
\end{eqnarray*}
The third equality owes to Lemma \ref{omega}.
\end{proof}

In order to prove the direction $(2) \Rightarrow (1)$ of 
Theorem \ref{ame},
we shall prove a more general statement as follows.
\begin{th}\label{AM}
Let $A$ be a finite von Neumann algebra with a fixed trace $\tr,$ and
$M \subset A$ be a hyperfinite II$_1$ factor.  Then we have
$[A^{\omega} : M^{\omega}] = [M' \cap A^{\omega} : M' \cap M^{\omega}].$
\end{th}

First we consider the case when $[A : M] < \infty,$
and use the idea of the Pimsner--Popa basis.
Here we recall their statements as below.

\begin{prop}\label{PPbasis}{\bf \cite[Proposition 1.3]{PP}}
Let ${\cal N} \stackrel{F}{\subset} {\cal M}$ be 
II$_1$ factors with the trace-preserving conditional expectation
$F$ from ${\cal M}$ onto ${\cal N}.$  When $[{\cal M} : {\cal N}] 
< \infty,$ there exists
a family $\{ m_j \}_{1 \leq j \leq n+1}$ of elements in ${\cal M},$ where
$n$ is the integer part of $[{\cal M} : {\cal N}],$
such that
 
(a)$F(m_j^* m_k) = 0, \quad j \neq k,$ 

(b)$F(m_j^* m_j) = 1, \quad 1 \leq j \leq n,$

(c)$F(m_{n+1}^* m_{n+1})$ is a projection of trace $[{\cal M} :
 {\cal N}] - n.$
\end{prop}

Now we fix some notations as follows.
Since $M$ is hyperfinite we may 
represent $M = \otimes_{n=1}^{\infty} M_2({\bf C})$ and set 
$A_m := \otimes_{n=1}^m M_2({\bf C}),
A := (M' \cap M_{\infty})' \cap M_{\infty}$ and $B := J_A M' J_A$ on $L^2(A).$
We easily notice that the following squares are commuting squares.
$$
\begin{array}{ccc}
M^{\omega} & \subset & A^{\omega} \\
\cup &       & \cup \\
\vdots & & \vdots \\
\cup & & \cup \\
A_k' \cap M^{\omega} & \subset & A_k' \cap A^{\omega} \\
\cup &       & \cup \\
A_{k+1}' \cap M^{\omega} & \subset & A_{k+1}' \cap A^{\omega} \\
\cup &       & \cup \\
\vdots & & \vdots \\
\cup && \cup \\
M' \cap M^{\omega} & \subset & M' \cap A^{\omega} 
\end{array}
$$
Our aim is to show the non-degeneracy of the next commuting square.
$$
\begin{array}{ccc}
M^{\omega} & \subset & A^{\omega} \\
\cup &       & \cup \\
M' \cap M^{\omega} & \subset & M' \cap A^{\omega}. 
\end{array}
$$
One can find a similar situation in 
\cite[Proposition 2.2 and Theorem 2.9]{P},
but we do not assume the factoriality of $A$ nor the finiteness
of the index.
The hyperfiniteness of $M$ is essential in main Theorem \ref{ame}.
We use the above presentation of $M$ throughout the proof.

\begin{proof}{\bf (Theorem \ref{AM}, the finite index case)} 
Since $M$ is a II$_1$ factor   % we have \begin{eqnarray*}
%{\bf C} = M' \cap M  &=& 
%((A_k' \cap M) \otimes A_k)' \cap ((A_k' \cap M) \otimes A_k) \\
%&=& \{(A_k' \cap M)' \cap (A_k' \cap M) \} \otimes (A_k' \cap A_k). 
%\end{eqnarray*} Then 
and $A_k \simeq M_{2^k}({\bf C}),$
the relative commutant $A_k' \cap M$ is also a II$_1$ factor.
Since $\dim (A' \cap A) \leq \dim (M' \cap A) < \infty,$
there exist a finite number of minimal central projections 
$q_1, \cdots, q_a \in A$
such that $A = Aq_1 \oplus \cdots \oplus Aq_a.$
We set each factor $A^{(i)} := Aq_i$ for $1 \leq i \leq a$ and 
$n_i$ to be the integer part of $[A^{(i)}: Mq_i].$
Since
\begin{eqnarray*} 
[(A_kq_i)' \cap A^{(i)}: (A_kq_i)' \cap Mq_i ] &\leq &
[A^{(i)} : (A_kq_i)' \cap Mq_i] \\
&=& [A^{(i)} : Mq_i] \times [Mq_i : (A_kq_i)' \cap Mq_i] \\
&\leq & [A : M] \times 2^{2k} < \infty,
\end{eqnarray*}
Proposition \ref{PPbasis} implies the existence of an orthonormal basis 
$\{m_{i,j}^{(k)}\}_{0 \leq j \leq n_i}$ of
$(A_kq_i)' \cap Mq_i \subset (A_kq_i)' \cap A^{(i)}$
such that

(1) $E_{Mq_i}(m_{i,j}^{(k) \, *} m_{i,j'}^{(k)}) = 0, \quad j \neq j',$

(2) $E_{Mq_i}(m_{i,j}^{(k) \, *} m_{i,j}^{(k)}) = q_i, \quad 
0 \leq j < n_i,$

(3) $E_{Mq_i}(m_{i,n_i}^{(k) \, *} m_{i,n_i}^{(k)})$ is a projection
of trace $[A^{(i)} : Mq_i] - n_i,$

(4) $\sum_{j=0}^{n_i} m_{i,j}^{(k)} m_{i,j}^{(k) \, *} = 
[A^{(i)} : Mq_i]q_i.$

Because 
$$A^{(i)} \simeq A_kq_i \otimes ((A_kq_i)' \cap A^{(i)}) \subset
\overline{\rm sp} \; Mq_i((A_kq_i)' \cap A^{(i)}),$$
the square
$$\begin{array}{ccc}
Mq_i & \stackrel{E|_{A^{(i)}}}{\subset} & A^{(i)} \\
\cup &                             & \cup \\
(A_kq_i)' \cap Mq_i & \subset & (A_kq_i)' \cap A^{(i)}
\end{array}$$
is a non-degenerate commuting square. 
Therefore, $ \{m_{i,j}^{(k)}\}_j$
is also an orthonormal basis of $Mq_i \subset A^{(i)},$
and $[A^{(i)} : Mq_i] = [(A_kq_i)' \cap A^{(i)} : (A_kq_i)' \cap Mq_i].$

We set $m_{i,j} := \{ m_{i,j}^{(k)} \}_k,$
then $\{ m_{i,j} \}_j$ is an orthonormal basis of both
$(Mq_i)' \cap M^{\omega}q_i \subset (Mq_i)' \cap A^{(i) \, \omega}$
and $M^{\omega}q_i \subset A^{(i) \, \omega}.$

Therefore, we have
\begin{eqnarray*}
[M' \cap A^{\omega} : M' \cap M^{\omega}]
&=& 
\sum_{i=1}^a [(Mq_i)' \cap A^{(i) \, \omega} : (Mq_i)' \cap M^{\omega}q_i] \\
&=& \sum_{i=1}^a \frac{ \tr (\sum_{j=0}^{n_i} m_{i,j}m_{i,j}^*)}{ \tr (q_i)}
= \sum_{i=0}^a [A^{(i) \, \omega} : M^{\omega}q_i] 
= [A^{\omega} : M^{\omega}].
\end{eqnarray*}
\end{proof}

\begin{proof}{\bf (Theorem \ref{AM}, the infinite index case with
$\dim (M' \cap A) = \infty$)} 
Since $A$ is a finite von Neumann algebra and 
$\dim (M' \cap A) = \infty,$
for any $\e > 0,$ there exists a non-zero projection $e \in M' \cap A$
such that $0 \neq \tr (e) \leq \e.$
Since the square$$
\begin{array}{ccc}
M & \subset & A \\
\cup &     &  \cup \\
M'\cap M = {\bf C} & \subset & M' \cap A
\end{array} $$
is a commuting square, we have
$0 \neq E_M (e) = \tr (e) \leq \e 1_M.$
Then if we set $\tilde{e} := \{e\}_k \in M' \cap A^{\omega},$
we have $\tilde{e} \neq 0$ and 
$E_{M' \cap M^{\omega}} (\tilde{e}) = E_{M^{\omega}}(\tilde{e})
\leq \e.$
By Proposition \ref{infty}, we obtain
$[M' \cap A^{\omega} : M' \cap M^{\omega}]_{PP} = \infty.$
By Proposition \ref{dim} (2), we have
$[M' \cap A^{\omega} : M' \cap M^{\omega}] = \infty.$
\end{proof}

%\begin{prop}\label{triv rel comm}{\bf \cite[pages 71, 72]{PP}}
%When $[M : A] = \infty$ and $M' \cap A = {\bf C},$ in particular,
%$A$ is a factor, then, 
%for any $\e >0,$ there exists a non-zero projection $e \in A$
%such that $E_M (e) \leq \e 1_M.$ \end{prop}
%We use this useful proposition without proof (see \cite[page 71, 72]{PP}).
%Later,we reduce the case of $[M : A] = \infty$ and $M' \cap A \neq {\bf C}$
%to Proposition \ref{triv rel comm}.
The next proposition will play an important role for the rest of the
proof
of Theorem \ref{AM}.  Pimsner and Popa have mentioned this statement as 
a Remark to Theorem 2.2 and Lemma 2.3 in \cite{PP} and say that
this follows from a maximality argument.
We shall include a full proof for the sake of completeness.
\begin{prop} \label{maximal}{\bf \cite[Remark 2.4]{PP}}
If $[A : M] = \infty$ and $M' \cap A = {\bf C},$ in particular,
$A$ is a factor, then,  
for any $\e \in (0, \, 1),$ there exists a non-zero projection $e \in A$
such that 
$$\tr (\chi_{\{ \e \}}(E_M (e))) \geq 1- \e.$$
\end{prop}

We prove that the maximal projection which satisfies $E_M (e) \leq \e 1_M$
is a one we desire.
This is based on the techniques in \cite[Theorem 2.2]{PP}.

\begin{proof}
By Zorn's lemma,
there exists a maximal family of mutually orthogonal projections 
$\{ e_i \}_{i \in I}$
in $A$ such that $E_M (e) \leq \e 1_M$
where we have set $0 \neq e := \sum_{i \in I} e_i.$

Suppose that the conclusion of the proposition dose not hold, i.e.,
$$\tr (\int_0^{\e} \chi_{[0, \, \e)} (\lambda) dE(\lambda)) > \e,$$
where $\int_0^{\e} \lambda dE( \lambda )$ is the spectral decomposition
of $E_M (e).$
Then there exists $\e_0 \in (0, \, \e)$ such that 
$$\tr (\int_0^{\e} \chi_{[0, \, \e_0 ]} (\lambda) dE(\lambda)) > \e.$$
We set $f_0 := \int_0^{\e} \chi_{[0, \, \e_0 ]} (\lambda) dE(\lambda) \in M.$
We need the following claim.

\begin{cl}\label{chi}
There exists a non-zero projection $q \in A$ such that
$q \leq 1-e,$ 
$q \leq f_0$ and
$E_M (q) \leq (\e - \e_0) f_0.$
\end{cl}

If we accept this claim we easily obtain a proof of this proposition as
 below.  We have
\begin{eqnarray*}
E_M (q \vee e) &=& E_M (q) +E_M (e) 
\leq (\e - \e_0)f_0 + E_M (e) \\ &\leq &
(\e - \e_0)f_0 + E_M (e)(1-f_0) + E_M (e)f_0 \\
&\leq & (\e - \e_0)f_0 + \e (1-f_0) + \e_0 f_0 = \e, 
\end{eqnarray*}  which contradicts the maximality of $\{ e_i \}_{i \in I}.$
\end{proof}

\begin{proof}{\bf (Claim \ref{chi})}
We set $q_0 := f_0 \wedge (1-e) \in A,$ then we obtain
\begin{eqnarray*}
\tr (q_0) &=& \tr (f_0) + \tr (1-e) - \tr (f_0 \vee (1-e)) \\
&\geq & \tr (f_0) + \tr (1-e) - 1 >
\e +1- \e -1 =0,
\end{eqnarray*}
which means $q_0 \neq 0.$
We denote $J_A M' J_A$ on $L^2(A)$ by $B$ 
and the Jones projection in $B$ by $e_0,$
then $B$ is a type II$_{\infty}$ factor.
Let $\phi$ be a semifinite trace on $B$ satisfying $\phi (e_0) = 1.$
Since $0 < \tr (q_0) \leq 1$ and 
$(q_0 A q_0)' \cap (q_0 B q_0) = q_0(A' \cap B)q_0 = 
{\bf C}q_0 \simeq {\bf C}$ by
 \cite[Lemma 2.1]{P2}, we know that $q_0 B q_0$ is a
type II$_{\infty}$ factor. 
Since $\tr (E_M (q_0)) = \tr (q_0) \neq 0,$ we have 
$(e_0 q_0) (e_0 q_0)^* = e_0 q_0 e_0 =
E_M (q_0)e_0 \neq 0.$ Therefore,
$q_0 e_0 q_0 = (e_0 q_0)^* (e_0 q_0) \neq 0$ and 
$q_0 e_0 q_0 \in (q_0Bq_0)_+.$ 
By the normalization of $\phi,$ we have
\begin{eqnarray*}
0 &<& \Vert q_0 e_0q_0 \Vert_{\phi} = \phi (q_0e_0q_0e_0q_0)\\
&=& \phi (e_0q_0e_0q_0e_0) = \phi (E_M(q_0)e_0E_M(q_0)e_0) \\
&=& \phi (E_M(q_0)^2 e_0) = \tr (E_M(q_0)^2) < \infty.
\end{eqnarray*}

Next we apply \cite[Lemma 2.3]{PP} to the inclusion
$q_0 A q_0 \subset q_0 B q_0$ with 
$\e$ replaced by 
$\sqrt{\tr (q_0)} (\e - \e_0)^2 / \Vert q_0 e_0 q_0 \Vert_{\phi}.$
(We repeat exactly the same arguments as in \cite[page 72]{PP} below.)
We obtain projections $f_1, \cdots, f_n \in q_0 A q_0$ such that
$\sum_{i=1}^n f_i = q_0$ and
\begin{eqnarray*}
\Vert \sum_{i=1}^n f_i (q_0 e_0 q_0) f_i \Vert_{\phi}^2
&< & \frac{\tr (q_0) (\e - \e_0)^4}{\Vert q_0 e_0 q_0 \Vert_{\phi}^2}
 \Vert q_0 e_0 q_0 \Vert_{\phi}^2    \\ 
&=& \tr (q_0) (\e - \e_0)^4 = \sum_{i=1}^n \tr (f_i)(\e - \e_0)^4.
\end{eqnarray*}
Since the left hand side equals to 
$\sum_{i=1}^n \Vert f_i q_0 e_0 q_0 f_i \Vert_{\phi}^2 
= \sum_{i=1}^n \Vert f_i e_0 f_i \Vert_{\phi}^2,$
there exists $j$ such that 
$$\Vert f_j e_0 f_j \Vert_{\phi}^2 < \tr (f_j) (\e - \e_0)^4.$$
We set $p := \chi_{[0,\, (\e - \e_0)]} (E_M (f_j)) \in M.$
By the normalization of $\phi$ and $0 < \e - \e_0 < 1,$ we have
\begin{eqnarray*}
\tr (f_j) (\e - \e_0)^4 &>& \Vert f_je_0f_j \Vert_{\phi}^2
= \phi (f_je_0f_je_0f_j) \\
&=& \phi (e_0f_je_0f_je_0) = \Vert e_0f_je_0 \Vert_{\phi}^2  \\
&=& \Vert E_M (f_j)e_0 \Vert_{\phi}^2 = \Vert E_M (f_j) \Vert_2^2 \\
&=& \tr (E_M (f_j)E_M (f_j)^*)
\geq \tr ((1-p)E_M (f_j)E_M (f_j)^*) \\
&\geq & (\e - \e_0)^2 \tr (1-p) 
\geq (\e - \e_0)^4 \tr (1-p), 
\end{eqnarray*}
which means $\tr (p+f_j) > 1.$
If we set $q := p \wedge f_j \in A,$ 
we have $$\tr (q) = \tr (p) + \tr (f_j) - \tr (p \vee f_j)
\leq \tr (p) + \tr (f_j) -1 > 0.$$
Thus,
$$0 \neq q = p \wedge f_j \leq f_j \leq q_0 \leq 1-e,$$
$$0 \neq q \leq q_0 \leq f_0$$ and 
\begin{eqnarray*}
E_M (q) &=& E_M (p \wedge f_j) \leq E_M (pf_jp) \\
&=& pE_M(f_j)p \leq (\e - \e_0)p.
\end{eqnarray*}
Multiplying the both hand sides of the third inequality by $f_0$ 
from the left and the right, we obtain
$$E_M (q) = E_M (f_0 q f_0) \leq (\e - \e_0)f_0 p f_0
\leq (\e - \e_0)f_0.$$
\end{proof}

We need the next %two 
lemma as a preparation for Lemma \ref{kityaku}.
%\begin{lm}\label{JPnotameni}
%Let $Q$ be a II$_1$ factor and $p_1, \cdots, p_n \in Q$
%be projections such that $\sum_{i=1}^n p_i = 1_Q.$
%Then we have 
%$$2^{n-1}   \sum_{i=1}^n q_ixq_i \geq x.$$
%\end{lm}
%
%\begin{proof}
%We prove it by induction on $k.$
%When $k=2,$ we have
%$$(p_1-p_2)x(p_1-p_2) \geq 0,$$ i.e.,
%$$p_1xp_1+p_2xp_2 \geq p_1xp_2+p_2xp_1,$$ thus
%\begin{eqnarray*}
%2(p_1xp_1+p_2xp_2) &\geq & p_1xp_1+p_2xp_2+p_1xp_2+p_2xp_1\\
%&=& (p_1+p_2)x(p_1+p_2)=x.
%\end{eqnarray*}
%We suppose the inequality holds when $n=k.$
%Since $(p_k-p_{k+1})x(p_k-p_{k+1}) \geq 0,$
%we have 
%$$2(p_kxp_k+p_{k+1}xp_{k+1}) \geq (p_k+p_{k+1})x(p_k+p_{k+1}).$$
%By the assumption of induction, we have
%\begin{eqnarray*}
%x & \leq & 2^{k-1}(\sum_{i=1}^{k-1}p_ixp_i +(p_k+p_{k+1})x(p_k+p_{k+1})) \\
%& \leq & 2^{k-1}(\sum_{i=1}^{k-1}p_ixp_i +2(p_kxp_k+p_{k+1}xp_{k+1})) \\
%& \leq & 2^k (\sum_{i=1}^k p_ixp_i) 
%\end{eqnarray*}
%\end{proof}

\begin{lm}\label{JP}
Let $P \subset Q$ be type II$_1$ factors with
$\dim (P' \cap Q) < \infty.$
Let $p_1, \cdots, p_n \in P' \cap Q$ be projections
such that $\sum_{i=1}^n p_i = 1.$
Then we obtain the following identity.
$$[Q : P] = \sum_{i=1}^n [p_iQp_i : Pp_i]/ \tr_Q (p_i)$$
\end{lm}

\begin{proof}
When $[Q : P] < \infty,$ 
Jones has already proved in \cite[Lemma 2.2.2]{J2}.

When $[Q : P] = \infty,$ we include a proof here
for the sake of completeness,
though it has been noted in \cite[page 61]{PP}.

For any $x \in Q,$ we have
\begin{eqnarray*}
E_P(p_ixp_i) &=& \sum_{j=1}^n p_jE_P(p_ixp_i)p_j
\geq p_iE_P(p_ixp_i)p_i \\
&=& E_{Pp_i}(p_ixp_i) \tr (p_i)
\geq \frac{\tr (p_i) p_ixp_i}{[p_iQp_i : Pp_i]_{PP}}.
\end{eqnarray*}
For any $x \in Q$ and $y \in P,$ we have
$$\tr_P(E_P(p_ixp_j)y) = \tr_Q (p_ixp_jy) = \tr_Q(xyp_jp_i).$$
Therefore, for any $x \in Q,$ we obtain
$$E_P(x)= \sum_{i=1}^n E_P(p_ixp_i).$$
By %Lemma \ref{JPnotameni} and 
the above arguments, 
if $[p_iQp_i : Pp_i]_{PP} = [p_iQp_i : Pp_i] < \infty$ for all $i,$
we have
\begin{eqnarray*}
E_P(x) &=& \sum_{i=1}^n E_P(p_ixp_i)
       \geq \sum_{i=1}^n \frac{\tr (p_i) p_ixp_i}{[p_iQp_i : Pp_i]} \\
     & \geq & \min_{1 \leq i \leq n} \left( \frac{\tr (p_i)}{[p_iQp_i : Pp_i]}
\right) \sum_{i=1}^n p_ixp_i
%       \geq 2^{-n+1} 
\geq \frac{1}{n}
\min_{1 \leq i \leq n} \left( \frac{\tr (p_i)}{[p_iQp_i : Pp_i]} \right) x,
\end{eqnarray*}
which contradicts $[Q : P] = \infty.$
The last inequality holds by
\cite[Definition 3.7.5]{GHJ} and \cite[Corollaire 2.3]{Jol}.
Thus, there exists an index $i_0$ such that
$[p_{i_0}Qp_{i_0} : Pp_{i_0}] = \infty.$
\end{proof}

Thanks to the following lemma, we can reduce the case of $\dim (M' \cap
A) < \infty$ to the case of $\dim (M' \cap
A) =1,$ which is the assumption of Proposition \ref{maximal}. 

\begin{lm}\label{kityaku}
Let $A$ be a finite von Neumann algebra with a fixed trace $\tr,$ and
$M \subset A$ be a type II$_1$ factor with 
$[A : M] = \infty$ and $\dim (M' \cap A) < \infty$.  Then there exist
a minimal central projection $q \in A' \cap A$ and 
a minimal projection $p \in (Mq)' \cap (qAq)$ such that
$qAq$ is a type II$_1$ factor,
$[pqAqp : pqMqp] = \infty$ and 
$(pqMqp)' \cap (pqAqp) = {\bf C}.$
\end{lm}
\begin{proof}
Since $A$ is a finite von Neumann algebra and 
$\dim (A' \cap A) \leq \dim (M' \cap A) < \infty,$
there exist minimal central projections $q_1, \cdots, q_n \in A' \cap A$
such that $A = \oplus_{i=1}^n q_i A q_i.$ 
We have
${}_M L^2(A) \simeq \oplus_{i=1}^n {}_M L^2 (q_i A q_i)$
as a left $M$-module, that is,
\begin{eqnarray*}
\infty &=& [A : M]  
= \sum_{i=1}^n \dim_M L^2 (q_i A q_i) \\
&=& \sum_{i=1}^n \dim_{Mq_i} L^2 (q_i A q_i) 
= \sum_{i=1}^n [q_i A q_i : Mq_i].
\end{eqnarray*}
Then there exists an index $j$ such that 
$[q_j A q_j : Mq_j] = \infty.$
Since we have 
\begin{eqnarray*}
\dim ((Mq_j)' \cap (q_jAq_j)) &=& \dim (M' \cap (q_jAq_j)) 
\leq \dim (\oplus_{i=1}^n (M' \cap q_iAq_i)) \\
&=& \dim(M' \cap A) < \infty,
\end{eqnarray*}
there exists a finite number of minimal projections
$p_1, \cdots , p_m \in (Mq_j)' \cap (q_j A q_j)$ such that
$\sum_{i=1}^m p_i = 1.$
By Lemma \ref{JP}, we obtain
$$\infty = [q_jAq_j : Mq_j] =
\sum_{i=1}^m[p_iq_jAq_jp_i : Mq_jp_i]/ \tr (p_i).$$
Then there exists an index $i$ such that $[p_iq_jAq_jp_i : Mq_jp_i] = \infty$ 
and $(p_iq_jMq_jp_i)' \cap (p_iq_jAq_jp_i) = {\bf C}$ by \cite[Lemma
 2.1]{P2} and the minimality of $p_j.$
\end{proof}

\begin{lm} \label{condi}
Let ${\cal N}$ be a type II$_1$ factor and ${\cal M} \supset {\cal N}$ be a
finite von Neumann algebra.  If $f_0 \in {\cal N}' \cap {\cal M}$
is a projection and $f \leq f_0$ is a projection, then we obtain the
following.
$$E_{\cal N}(f)f_0 = \tr (f_0) E_{{\cal N}f_0}(f)$$
\end{lm}

This is also noted in \cite[page 71]{PP}, but we include a proof.

\begin{proof}
For any $x \in {\cal N}$ we have
\begin{eqnarray*}
\tr_{{\cal N}f_0}(E_{{\cal N}f_0}(f)x) \tr_{\cal M} (f_0)
&=&  \tr_{{\cal N} f_0}(ff_0x) \tr_{\cal M}(f_0) \\
&=&  \tr_{\cal M}(ff_0x) = \tr_{\cal M}(fx) \\
&=& \tr_{\cal N} (E_{\cal N}(f)x) = \tr_{{\cal N}f_0}(E_{\cal N}(f)xf_0).
\end{eqnarray*}
Since $f_0 \in {\cal N}' \cap {\cal M},$ we obtain the identity
$E_{{\cal N}}(f)f_0 = \tr (f_0) E_{{\cal N}f_0}(f).$
\end{proof}

\begin{proof}{\bf (Theorem \ref{AM}, the infinite index case with 
$\dim (M' \cap A) < \infty$)}
By Proposition \ref{kityaku},
there exist projections $p \in A' \cap A$ and
$q \in (Mq)' \cap (qAq)$ such that
$[pqAqp : pqMqp] = \infty$ and $(pqMqp)' \cap (pqAqp) = {\bf C}.$
And by Lemma \ref{maximal},
for any $\e >0,$ there exists a non-zero projection
$e \in pqAqp$ such that
$E_{pqMqp}(e) \leq \e qp$ and
$\tr_{pqMqp}(\chi_{\{ \e \}}(E_{pqMqp}(e))) \geq 1-\e .$
Then by Lemma \ref{condi}
(here we use the idea of \cite{PP}), we have
$E_ M(e)qp = \tr_A(qp)E_{pqMqp}(e) \leq \e \tr_A (qp)qp.$ 

Since $qp \in M' \cap A,$ $M$ is isomorphic to $Mqp,$ thus
$E_M(e) \leq \e \tr_A (qp)1_M.$
We also have,
\begin{eqnarray*}
\tr_M(E_M(e))& =& \tr_A(e) = \tr_A(eqp) = \tr_{pqAqp}(eqp)\tr_A (pq) \\
&=& \tr_{pqMqp}(E_{pqMqp}(e)) \tr_A(pq)
\geq \e (1- \e )\tr_A (pq) \neq 0.
\end{eqnarray*}

Since $pq \in M' \cap A = (A_k' \cap M)' \cap (A_k' \cap A),$
we have 
$$(pq(A_k' \cap M)qp)' \cap (pq(A_k' \cap A)qp) = {\bf C}$$
and
$$[pq(A_k' \cap A)qp : pq(A_k' \cap M)qp] = \infty$$
in the same way.  Then 
for any $\e,$ there exists a non-zero projection $e_k \in pq(A_k' \cap A)qp 
\subset A_k' \cap A$
such that $E_{A_k' \cap M}(e_k) \leq \e \tr_A (qp) 1$ and
$\tr (e_k) \geq \e (1- \e ) \tr (pq) \neq 0.$

By setting $e := \{ e_k \}_k \in M' \cap A^{\omega},$
we have $e \neq 0$ and $E_{M' \cap M^{\omega}}(e) 
\leq \e 1_{M' \cap M^{\omega}}.$
Thanks to Propositions \ref{infty} and \ref{dim}, we obtain 
$[M' \cap A^{\omega} : M' \cap M^{\omega}] = \infty.$
\end{proof}

\begin{proof}{\bf (Theorem \ref{ame} (2) $\Rightarrow$ (1))}
Thanks to Theorem \ref{AM} and Lemma \ref{omeg}, 
we obtain 
$$(M' \cap M_{\infty})' \cap M_{\infty} = A = M.$$
\end{proof}
\begin{Remark}
So far we have considered only the larger algebra of
$M' \cap N^{\omega} \subset M' \cap M^{\omega}$
and proved the equivalence between the double commutant property of
$ M' \cap M^{\omega}$ and Popa's strong amenability.
It is known that when $N \subset M$ has an infinite depth, we have
$[M' \cap M^{\omega} : M' \cap N^{\omega}] = \infty.$
Thus if we define $P_2$ by the Jones basic construction of 
$M' \cap N^{\omega} \subset M' \cap M^{\omega},$
it does not make sense to consider the double commutant properties of
$P_k$ $(k \geq 2),$ because $P_2 \not\subset M_{\infty}^{\omega}.$

As for the smaller algebra, we have not considered yet.
In general, the condition
$((M' \cap N^{\omega})' \cap M_{\infty}^{\omega})' \cap M_{\infty}^{\omega}
= M' \cap N^{\omega}$ does not imply the strong amenability.
(We recall that the converse direction always holds, see 
\cite[Section 15.5]{EK}.)
For example, 
let $R_0$ be a hyperfinite II$_1$ factor.  We set
$G := PSL(2, {\bf Z}) \simeq {\bf Z}/ 2{\bf Z} \ast {\bf Z}/3{\bf Z},$
$R := \otimes_{g \in G} R_0 ( \simeq R_0),$ and $\alpha$ 
to be an outer action of $G$ on $R$ defined by the Bernoulli shift.
We restrict the action $\alpha$ to ${\bf Z}/ 2{\bf Z}$ and 
${\bf Z}/3{\bf Z}$ regarded as subgroups of $G,$ and set 
$N := R^{{\bf Z}/ 2{\bf Z}}$ and
$M := R \rtimes ({\bf Z}/3{\bf Z}).$
Since $G$ is non-amenable, by \cite[Proposition 2]{Jremark},
we have $N^{\omega} \cap M' = (R_{\omega})^G = {\bf C}.$
This example in \cite[page 211]{B} is due to Jones and based on
\cite{Jremark}.
Then $N^{\omega} \cap M' = {\bf C} = (N^{\omega} \cap M')^{cc}.$
If the subfactor were strong amenable, the generating property
(see \cite[Theorem 4.2.1]{Pacta}) would imply
\begin{eqnarray*}
(N \subset M) & \simeq & 
\left( \overline{\bigcup_{k=1}^{\infty} N_k' \cap N} \subset 
\overline{\bigcup_{k=1}^{\infty} N_k' \cap M} \right) \\
&\simeq & \left( \overline{\bigcup_{k=1}^{\infty} (R \otimes N_k)' 
\cap (R \otimes N)} \subset 
\overline{\bigcup_{k=1}^{\infty} 
(R \otimes N_k)' \cap (R \otimes M)} \right) \\
&\simeq & (R \otimes N \subset R \otimes M),
\end{eqnarray*}
where $\cdots N_2 \subset N_1 \subset N \subset M$
is a generating tunnel.
Then we have
\begin{eqnarray*}
{\bf C} &\simeq & N^{\omega} \cap M' \simeq
(R \otimes N)^{\omega} \cap (R \otimes M)' \\
& \supset & R_{\omega} \otimes {\bf C} \simeq R_{\omega},
\end{eqnarray*}
which is a contradiction.
\end{Remark}


\begin{thebibliography}{99}
\bibitem[BDH]{BDH} M. Baillet, Y. Denizeau, and J.-F. Havet,
{\it Indice d'une esperence conditionelle}, 
Compos. Math. {\bf 66} (1988), 199--236.

\bibitem[B]{B} D. Bisch,
{\it A note on intermediate subfactors},
Pac. J. Math. {\bf 163} (1994), 201--216.

\bibitem[EK]{EK} D. E. Evans and Y. Kawahigashi,
{\it Quantum symmetries on operator algebras},
Oxford University Press, (1998).

\bibitem[GHJ]{GHJ}
F. Goodman, P. de la Harpe, \& V. F. R. Jones,
{\it Coxeter graphs and towers of algebras}, MSRI publications {\bf 14}
Springer, (1989).

\bibitem[H]{H} U. Haagerup,
{\it Operator valued weights in von Neumann algebras} II,
J. Funct. Anal. {\bf 33} (1979), 339--361.

\bibitem[J1]{Jremark} V. F. R. Jones,
{\it A converse to Ocneanu's theorem},
J. Operator Theory {\bf 10} (1983), 61--63.

\bibitem[J2]{J2} V. F. R. Jones, 
{\it Index for subfactors}, Invent. Math. {\bf 72} (1983), 1--25.

\bibitem[Jol]{Jol} P. Jolissaint,
{\it Indice d'esperances conditionnelles et algebres de von Neumann
 finies}, Math. Scand. {\bf 68} (1991), 221--246.

\bibitem[K]{Ko} H. Kosaki,
{\it Extension of Jones theory on index to arbitrary factors},
J. Funct. Anal. {\bf 66} (1986), 123--140.

\bibitem[L]{L}  R. Longo,
{\it Index of subfactors and statistics of quantum fields. I},
Commun. Math. Phys. {\bf 126} (1989), 217--247.

\bibitem[LR]{LR}  R. Longo and K.-H. Rehren,
{\it Nets for subfactors},
Rev. Math. Phys. {\bf 7} (1995), 567--597.

\bibitem[O1]{O1}  A. Ocneanu,
{\it Quantized group string algebras and Galois theory for algebras},
in ``Operator algebras and applications, Vol. 2
(Warwick, 1987),''
London Math. Soc. Lect. Note Series Vol. 136,
Cambridge University Press,
1988, pp. 119--172.

\bibitem[O2]{O2}  A. Ocneanu,
{\it Quantum symmetry, differential geometry of finite graphs and
classification of subfactors}, University of Tokyo Seminary Notes {\bf 45}
(Notes recorded by Y. Kawahigashi), (1991).

\bibitem[PP]{PP} M. Pimsner and S. Popa,
{\it Entropy and index for subfactors}, Ann. Scient. \'{E}c. Norm. Sup.
4$^e$ s\'{e}rie, {\bf 19} (1986), 57--106.

\bibitem[P1]{P2} S. Popa,
{\it On a problem of R. V. Kadison on maximal abelian $*$-subalgebras in
factors}, Invent. Math. {\bf 65} (1981), 269--281.

\bibitem[P2]{P0} S. Popa,
{\it Orthogonal pairs of $*$-algebras in finite von Neumann
algebras}, J. Operator Theory {\bf 9} (1983), 253--268.

\bibitem[P3]{Paction} S. Popa,
{\it Sousfacteurs, actions des groupes et cohomologie},
C. R. Acad. Sci. Paris  {\bf 309} S\'{e}rie I, (1989), 771--776.

\bibitem[P4]{Pacta} S. Popa,
{\it Classification of amenable subfactors of type II}, Acta Math. {\bf 172}
(1994), 163--255.

\bibitem[P5]{P} S. Popa,
{\it Classification of subfactors and their endomorphisms}, CBMS 
Regional Conference Series in Mathematics {\bf 86}, (1995).

\bibitem[S]{S} S. Str\u atil\u a,
{\it Modular theory in operator algebras}, Abacus Press--Editura Academiei,
London, Bucharest, (1981).
\end{thebibliography}
\end{document}